\documentclass[12pt]{amsart}

\usepackage{CJK}

\usepackage{float}
\usepackage{subfig}
\usepackage{indentfirst}
\usepackage{mathrsfs}
\usepackage[top=1.2in,bottom=1.15in,left=1.45in,right=1.45in]{geometry}

\usepackage{amsfonts,latexsym,rawfonts,amsmath,amssymb,amsthm,times,a4wide}
\usepackage{graphicx}
\usepackage[plainpages=false]{hyperref}

\usepackage{enumerate}
\usepackage{bbm}

\numberwithin{equation}{section}
\DeclareMathOperator{\arccosh}{arccosh}

\usepackage{amscd}

\usepackage{eufrak}
\usepackage{euscript}
\usepackage{epsfig}

\usepackage{array}

\usepackage{color}
\usepackage{wasysym}
\usepackage{pdfsync}
\usepackage{stmaryrd}

\newcommand{\beq}{\begin{equation}}
\newcommand{\eeq}{\end{equation}}
\newcommand{\beqs}{\begin{eqnarray*}}
\newcommand{\eeqs}{\end{eqnarray*}}
\newcommand{\beqn}{\begin{eqnarray}}
\newcommand{\eeqn}{\end{eqnarray}}
\newcommand{\beqa}{\begin{array}}
\newcommand{\eeqa}{\end{array}}

\newtheorem{prop}{Proposition}[section]
\newtheorem{theo}[prop]{Theorem}
\newtheorem{lem}[prop]{Lemma}

\newtheorem{cor}[prop]{Corollary}
\newtheorem{rem}[prop]{Remark}

\newtheorem{defi}[prop]{Definition}

\newcommand{\R}{{\mathbb R}}

\newcommand{\T}{\mathcal{T}}

\newcommand{\IH}{\mathbb{H}}
\newcommand{\IK}{\mathbb{K}}
\newcommand{\CB}{{\mathcal B}}
\newcommand{\CL}{{\mathcal L}}
\newcommand{\tm}{\begin{theo}}
\newcommand{\tmd}{\end{theo}}
\newcommand{\co}{\begin{cor}}
\newcommand{\cod}{\end{cor}}
\newcommand{\prp}{\begin{prop}}
\newcommand{\prpd}{\end{prop}}

\newcommand \di{\,\mathrm{d}}

\begin{document}

\title{Hyperbolization and geometric decomposition of a class of 3-manifolds}
\author{Ke Feng, Huabin Ge, Yunpeng Meng}
\address{}
\curraddr{}
\email{}
\thanks{}
\date{\today}
\dedicatory{}
\begin{abstract}
Thurston's triangulation conjecture asserts that every hyperbolic 3-manifold admits a geometric triangulation into hyper-ideal hyperbolic tetrahedra. So far, this conjecture had only been proven for a few special 3-manifolds. In this article, we confirm this conjecture for a class of 3-manifolds. To be precise, let $M$ be an oriented compact 3-manifold with boundary, no component of which is a $2$-sphere, and $\T$ is an ideal triangulation of $M$. If $\T$ satisfies properly gluing condition, and the valence is at least 6 at each ideal edge and 11 at each hyper-ideal edge, then $M$ admits an unique complete hyperbolic metric with totally geodesic boundary, so that $\T$ is isotopic to a geometric ideal triangulation of $M$.

We use analytical tools such as combinatorial Ricci flow (CRF, abbr.) to derive the conclusions. There are intrinsic difficulties in dealing with CRF. First, the CRF may collapse in a finite time, second, most of the smooth curvature flow methods are no longer applicable since there is no local coordinates in $\T$, and third, the evolution of CRF is affected by certain combinatorial obstacles in addition to topology. To this end, we introduce the ideas as ``extending CRF", ``tetrahedral comparison principles", and ``control CRF with edge valence" to solve the above difficulties. In addition, the presence of torus boundary adds substantial difficulties in this article, which we have solved by introducing the properly gluing conditions on $\T$ and reducing the ECRF to a flow relatively easy to handle.
\end{abstract}
\maketitle

\section{Introduction}
\label{sec-introduction}
It is conjectured that every hyperbolic 3-manifold, with both cusps and totally geodesic boundary, admits a geometric decomposition into hyper-ideal tetrahedra. It is originated from the creative idea of constructing cusped 3-manifolds by pasting isometrically ideal hyperbolic tetrahedra \cite{CDGW2016}. On the contrary, Thurston thought that each cusped 3-manifold can be decomposed into ideal hyperbolic tetrahedra and he once used it as a correct premise when he proved hyperbolic Dehn filling theorem. Due to this reason, the conjecture is referred to as Thurston's triangulation conjecture in our paper. To some extent, the triangulation conjecture is one of the most unsolved problems in the field of 3D geometry and topology, and is proven only for certain special 3-manifolds. ``It is a difficult problem in general", as was pointed by
Gu\'{e}ritaud-Schleimer \cite{GS2010}, ``general results are known only when $M$ is restricted to belong to certain classes of manifolds: punctured-torus bundles, two-bridge link complements, certain arborescent link complements and related objects, or covers of any of these spaces" (see \cite{Aki}-\cite{Aki-SWY-2}, \cite{Gue-1,Gue-2,Ham-P,Jorgen,Lackenby,Nimer} for example). In this paper, we confirm this conjecture for a class of 3-manifolds.

\begin{theo}\label{mainthm}
Let $M$ be an oriented compact 3-manifold with boundary, no component of which is a $2$-sphere, and $\T$ is a properly glued ideal triangulation of $M$. Assume the valence of edges respectively satisfies $d_e \ge 6$ for each ideal edge and $d_e \ge 11$  for each hyper-ideal edge, then $M-\partial_t$ ($\partial_t$ is the total boundary components of $M$) admits a finite-volume complete hyperbolic structure with totally geodesic boundary, and $\T$ is isotopic to a geometric decomposition of $M$.
\end{theo}

By definition, the \emph{valence} of an edge $e$ is the number of tetrahedra pasted around $e$. A type $3$-$1$ ($4$-$0$, resp.) tetrahedron is a hyperbolic tetrahedron that has 1 (0, resp.) ideal vertex and 3 (4, resp.) truncated hyper-ideal vertices. The side connecting at least one ideal vertex is called an \emph{ideal edge}, and the side connecting two different hyper-ideal vertices is called a \emph{hyper-ideal edge}. In addition, the triangulation $\T$ is called \emph{properly glued}, if each of the tetrahedra in it is of type $3$-$1$ or $4$-$0$, and if two $3$-$1$ tetrahedra $\sigma_{1234}$ and $\sigma_{1'2'3'4'}$ (with ideal vertex $1$ or $1'$ resp.) are glued along the face $123$ and $1'2'3'$ ($k$ and $k'$ are glued for $1\leq k\leq 3$), then the proper gluing requires sides $24$ and $2'4'$ to be glued, $34$ and $3'4'$ to be glued. See Section \ref{sec-prop-triangulation} for details.

\begin{rem}
For any torus boundary component $C$ of the manifold $M$, there is a triangulation $\T_C$ induced by the triangulation $\T$. The link at each ideal (hyper-ideal, resp.) point in $\T$ is a torus (surface with genus at least 2, resp.). If $d_e\geq6$ at each ideal edge $e$ in $\T$, the degree at each vertex is at least 6 in $\T_C$. By Euler's Polyhedron Formula, it's easy to show that the number has to be 6. This means that $d_e \ge 6$ in the above theorem is actually equivalent to $d_e = 6$.
\end{rem}

The use of \emph{edge valence} to obtain a hyperbolic structure on $(M,\T)$ is a pioneering discovery of Costantino-Frigerio-Martelli-Petronio in \cite{CFMP}, where they showed that if all edge valencies are at least 6, then $(M,\T)$ admits a hyperbolic structure. Moreover, they conjectured that the original triangulation $\T$ is also geometric. If $\partial M$ consists of surfaces of genus at least 2, \cite{Feng2022} confirmed their conjecture with all $d_e\geq10$. If $\partial M$ has at least a torus boundary component, their conjecture is confirmed in this article with $d_e\geq$ 6 or 11 at each ideal edge or hyper-ideal edge $e$ respectively. Similar phenomenon had been observed in other cases. For example, for a triangulated closed surface $(X,\T)$, if each vertex degree is at least 7, then $X$ admits a hyperbolic circle packing realizing $\T$, so that $X$ is hyperbolic and $\T$ is geometric \cite{GeLin}. For a triangulated closed 3-manifold $(M,\T)$, if each vertex degree is at least 23, then $M$ admits a (virtual) hyperbolic sphere packing realizing $\T$ \cite{GeHua}. 
Recall the \emph{vertex degree} here is the number of triangles or tetrahedra pasted around this vertex in dimension 2 or 3 respectively. By \cite{Lackenby2000,Lackenby-AGT} and Corollary 1.3 in \cite{Feng2023}, we have

\begin{cor}
\label{thm-converge-imply-topo}
Suppose $M$ is an oriented compact 3-manifold with boundary, no component of which is a 2-sphere. If $M$ admits a properly glued ideal triangulation $\T$, so that $d_e \ge 6$ for each ideal edge and $d_e \ge 11$  for each hyper-ideal edge, then $M$ is irreducible, atoroidal and not Seifert fibred. Moreover, $\T$ supports an angle structure and is strongly $1$-efficient and therefore $1$-efficient and $0$-efficient.
\end{cor}

The main tool to prove Theorem \ref{mainthm} is the \emph{combinatorial Ricci flow} (\emph{CRF}, abbr.). In fact, to obtain the results in \cite{Feng2022,GeHua,GeLin}, the combinatorial Ricci or Yamabe flow plays a crucial role and links the combinatorics of $\T$ to the geometry of $M$ in an ingenious way. CRF on surfaces were creatively introduced by Chow-Luo \cite{Chow-Luo} aiming to repove Koebe-Andreev-Thurston's circle packing theorem. Luo also initiated a program \cite{Luo2005} to study the geometry and topology of 3-manifolds using the CRF. Follow-up studies have not progressed too well, mainly due to three essential difficulties: first, the CRF may collapse in a finite time, second, the smooth Ricci flow techniques are no longer applicable since there is no local coordinates in $\T$, and third, the evolution of CRF is affected by certain combinatorial obstacles in addition to topology. By introducing the ideas as ``extending CRF", ``tetrahedral comparison principles", and ``control CRF with edge valence" in \cite{Feng2022,GeHua,GeJShen,GeLin}, these difficulties were solved. In addition, the presence of torus boundary adds substantial difficulties in this paper. To address this, we introduce properly gluing conditions on $\T$, perfectly reducing the CRF to one relatively easy to handle.

In the final part of the introduction, we outline the definition of the CRF, as well as the basic ideas and main steps to prove Theorem \ref{mainthm}. It is well known that hyperbolic 3-manifolds can be produced by gluing hyperbolic tetrahedra. For a compact 3-manifold $M$ with ideal triangulation $\T$, replace each topological tetrahedron in $\T$ with a hyperbolic tetrahedron, and paste these hyperbolic tetrahedra isometrically along the original triangulation, we get a hyperbolic polyhedral manifold $M$ with cone singularity on each edge of $\T$. The metric is non-singular and thus hyperbolic exactly when all the cone angles are equal to $2\pi$. At this time, $M$ is hyperbolic and $\T$ is geometric.  Traditionally, with the help of Thurston's gluing equations, a hyperbolic metric with all cone angles equal to $2\pi$ is obtained. In contrast, Luo's program is to evolve the hyperbolic polyhedron metric using CRF such that all cone angles tend to $2\pi$. We give slightly more specific explanations and descriptions below.

Let $M$ be a compact 3-manifold with boundary, no component of which is a $2$-sphere, and $\T$ is an ideal triangulation of $M$ with each tetrahedron combinatorially equivalent to a type $3$-$1$ or $4$-$0$ hyperbolic tetrahedron. By definition, a decorated $3$-$1$ type tetrahedron is a type $3$-$1$ hyperbolic tetrahedron which has a decoration, i.e. a horosphere centered at the ideal vertex (see \cite{Feng2023,FengGLiu} for more explanations). The sides lie entirely in $\partial M$ are not considered as edges in $\T$. Hence an edge in $\T$ is also called an internal edge, which is either an ideal edge or a hyper-ideal edge.

\begin{defi}[decorated metric]
\label{def:t1}
A decorated hyperbolic polyhedral metric on $(M, \T)$ is obtained by replacing each tetrahedra by a decorated $3$-$1$ type or $4$-$0$ type hyperbolic tetrahedron according its type and replacing the affine gluing homeomorphisms by isometries preserving the decoration, i.e. gluing decorated hyperbolic tetrahedra along codimension-$1$ faces. By the construction, such a metric is determined by the signed lengths of the edges in $\T$. We denote
$$l=(l(e_1), \dots, l(e_m))$$
by the decorated edge length vector, where $E = \{e_1, \cdots, e_m\}$ is the set of edges in $\T$. In fact, we usually refer to $l$ as a decorated hyperbolic polyhedral metric, or \textbf{decorated metric} for short.
\end{defi}

Denote $\mathcal{L}(M, \mathcal{T})\subset \R^E$ by the set of all decorated metrics on $(M, \mathcal{T})$ parametrized by the decorated edge length vector $l$. Just like each metric tensor determines a curvature tensor in a Riemannian manifold, each decorated metric $l$ determines a combinatorial curvature as follows.

\begin{defi}[combinatorial Ricci curvature]
\label{def:curvature}
Let $l\in(M, \mathcal{T})$ be a decorated metric. The \textbf{combinatorial Ricci curvature} at an edge $e$ is defined as
$$K_e(l)=2\pi-C_e,$$
where $C_e$ is the cone angle at $e$, i.e. the sum of the dihedral angles around $e$ in tetrahedra incident to $e$. This provides the combinatorial Ricci curvature vector $K=K(l)=(K_{e_1}(l), \cdots, K_{e_m}(l))$.
\end{defi}

A decorated metric $l$ with zero curvature, if exists, is geometrically particular meaningful. It makes $M$ hyperbolic and meanwhile $\T$ geometric. Such metric is unique up to hyperbolic isometry and change of decorations. More generally, $l$ is uniquely determined by $K$ up to isometry and change of decorations. This phenomenon is known as rigidity and is demonstrated in \cite{Feng2023} for general type-mixed triangulations. Rigidity of type $4$-$0$, $0$-$4$ and $1$-$3$ triangulations had been derived in \cite{Luo2018} and \cite{FengGLiu} before. To find such a decorated metric $l$ with zero curvature, the CRF provides an extremely suitable tool.

\begin{defi}[CRF]
The \textbf{combinatorial Ricci flow} on $(M,\T)$ for decorated metrics $l(t)$ is a system of ODEs with the time parameter $t$ so that
\begin{equation}
\label{eq:luoflow}
\frac{d }{dt}l_i(t)=K_i(l(t)),\quad i\in E, \;\; t\geq 0.
\end{equation}
\end{defi}

The solution $l(t)$ may attain $\partial\mathcal{L}(M, \mathcal{T})$ in finite time. To overcome this finite time collapsing phenomenon of the CRF, we introduced the extended CRF
following the extension idea pioneered by \cite{Bobenko2015} and \cite{Luo2011}:
\begin{defi}[ECRF]
The \textbf{extended combinatorial Ricci flow} on $(M,\T)$ is defined as
\begin{equation}
\label{eq:luoflow-ECRF}
\frac{d }{dt}l_i(t)=\widetilde{K_i}(l(t)),\quad i\in E, \;\; t\geq 0,
\end{equation}
where $\widetilde{K}$ is the generalized curvature, and $l(t)$ are generalized metrics. See section \ref{sec-ECRF}.
\end{defi}
By extending the so called \emph{co-volume function}, which is convex, and is the Legendre-Fenchel dual transformation of the volume function (the total volume of all hyperbolic tetrahedra), we rewrite (\ref{eq:luoflow-ECRF}) as a negative gradient flow of the $H$-function (which differs the extended co-volume function by a linear term), and established the following \emph{fundamental theorem of ECRF} in \cite{Feng2023} with the help of convex optimization methods and techniques.

\begin{theo}[Fundamental Theorem of ECRF, \cite{Feng2023}]
For given $(M,\T)$, the solution $\{l(t)\}$ to the ECRF is unique, and exists for all time $t$. Moreover, there is a hyperbolic structure on $M$ so that $\T$ is isotopic to a geometric triangulation if and only if $l(t)$ converges.
\end{theo}

If $M$ has no torus boundary, we \cite{Feng2022} developed tetrahedral comparison principles and introduced the combinatorial conditions, such as ``\emph{edge valence}", to control the ECRF and make it converge to a hyperbolic structure and geometric triangulation. If $M$ has at least one torus boundary, which is also the case in this paper, ideal vertices and ideal edges will appear in $\T$, this will cause additional difficulties, making it impossible for the method in \cite{Feng2022} to be directly generalized to this paper. By analyzing the structure of the ECRF and using the properly gluing assumption, we reduce it to the following flow, which is only defined on hyper-ideal edges.
\begin{defi}[the reduced ECRF]
Let $E=E^*\sqcup E^h$, where $E^*$ is the set of ideal edges and $E^h$ is the set of hyper-ideal edges in $\T$. Suppose $\T$ is properly glued, the reduced ECRF is define as
\begin{equation}
\label{eq:reduced-ECRF-introduction}
\dfrac{d }{d t}l^h(t)=\widetilde{K}(l(t)), \quad t\geq 0.
\end{equation}
Here $l=(l^*,l^h)$, and $l^*$ is the restriction of $l$ on $E^*$, $l^h$ is the restriction of $l$ on $E^h$.
\end{defi}

The reduced flow (\ref{eq:reduced-ECRF-introduction}) is more subtle than it may seem at first glance. Since $\widetilde{K}$ is a function of $l$, but not of $l^h$, the equation (\ref{eq:reduced-ECRF-introduction}) is not meaningful in general. This forced us to consider properly glued triangulations. For such a triangulation $\T$, each $l^h$
determines a unique generalized metric $l$ assuming each ideal vertex triangle is equilateral with length $1$. Hence $\widetilde{K}$ is actually a function of $l^h$ if $\T$ is properly glued (and further satisfying the assumption that all ideal vertex triangles are equilateral with length $1$).

We use the reduced ECRF (\ref{eq:reduced-ECRF-introduction}) to prove Theorem \ref{mainthm}, with proof strategies ketched as follows.

Through the meticulous study of the geometry of $3$-$1$ tetrahedra in subsection \ref{sec-31proper}, Proposition \ref{surjection} and Corollary \ref{prop-equa-length-1} prove that for any $l^h\in\mathbb{R}^{E^h}_{>0}$, there exists a unique generalized metric $l\in\mathbb{R}^E$ so that $l=(l^*,l^h)$ and each ideal vertex triangle in $\T$ is an Euclidean equilateral triangle with length $1$. Since $l=(l^*, l^h)$ is uniquely determined by its $l^h$ part, the reduced ECRF (\ref{eq:reduced-ECRF-introduction}) is well defined. Moreover, we may rewrite the extended co-volume function $cov$ and $H$ as $cov^h$ and $H^h$ since they are functions of $l^h\in\mathbb{R}^{E^h}_{>0}$. By proving $\nabla_{l^h} H^h=-\widetilde{K}$, we know that (\ref{eq:reduced-ECRF-introduction}) is a negative gradient flow, and hence $H^h(l^h(t))$ is decreasing along (\ref{eq:reduced-ECRF-introduction}). For each $e\in E^h$, Proposition \ref{prop:c0} in this article and Proposition \cite{Feng2022} show that the cone angles $C_e$ is uniformly small if $l(e)$ is close to zero, implying that the solution $l^h(t)$ to (\ref{eq:reduced-ECRF-introduction}) will never touch the zero boundary of $\mathbb{R}^{E^h}_{>0}$ and then exits for all time $t\in[0,+\infty)$.

Choosing a small initial metric $l^h_0\in (0,\arccosh 2)^{E^h}$, we may prove that the solution $l^h(t)$ is compactly supported in $(0,\infty)^{E^h}$, i.e. there exists a constant $c > 0$ so that $l^h(t)\in (c,\arccosh 2)^{E^h}$ for all $t\geq 0$. To get this uniform bound, firstly, Corollary \ref{lem:longest} establishes an estimate for the dihedral angle at the longest hyper-ideal edge in a generalized $3$-$1$ type tetrahedron, then Theorem \ref{LONG} derives an upper bound estimate for $l^h(t)$ under the combinatorial assumption $d\geq11$. Secondly, Theorem \ref{bd} derives a lower bound estimate based on the previous upper bound estimate and Proposition \ref{prop:c0}.

Since $l^h(t)$ is compactly supported in $\mathbb{R}^{E^h}_{>0}$, using the monotonicity of $H^h(l^h(t))$ along the reduced ECRF (\ref{eq:reduced-ECRF-introduction}), there is a subsequence $l^h(t_n)$ which converges to a limit metric $l^h_{\infty}$ as $n$ tends to $\infty$. By Proposition \ref{surjection}, a zero curvature decorated metric $l_{\infty}$ could be reconstructed from $l^h_{\infty}$. By Proposition \ref{realmetric}, $l_{\infty}$ is a real metric, and determines a hyperbolic structure on $M$ so that $\T$ is isotopic to a geometric ideal triangulation.

We remark that compared to \cite{Feng2022,Feng2023}, we do not know whether the solution $l^h(t)$ is unique, or whether the Lyapunov function $H^h$ is convex. Thus the convergence of $l^h(t)$ ($t\rightarrow\infty$) can not be derived generally. Moreover, most of our results are still valid if $(M,\T)$ is a triangulated compact pseudo 3-manifold. It is also worth mentioning that before the introduction of the CRF, there had been a lot of creative and profound work that linked the combinatorial of the triangulation with the geometry of $M$. For instance, we refer to Thurston's combinatorial obstructions (see section 13.7 in \cite{Thurston2022}) for the existence of circle packings, and Rivin's seminal work \cite{Rivin1994,Rivin2003}.

The paper is organized as follows. In Section \ref{sec-31type}, we recall some results related to the geometry of type $3$-$1$ hyperbolic tetrahedra.
In Section \ref{sec-proper}, we prove some new geometric properties for type $3$-$1$ and $4$-$0$ hyperbolic tetrahedra. We further give some explanation of how the conclusion changes as the conditions change. At the end of this section, we give a further discussion of the co-volume function.
In Section \ref{sec-erf}, for $(M, \T)$, we will introduce ``properly gluing conditions" and recall the extended combinatorial Ricci flow and corresponding conclusions. In Section \ref{sec-reduce-ECRF}, we propose a reduced extended combinatorial Ricci flow and derive corresponding conclusions. In Section \ref{sec-proofthm}, we prove Theorem \ref{mainthm}.

\bigskip
\noindent {\bf Acknowledgements.}
The first author is supported by NSFSC, no. 2023NSFSC1286.
The second author is supported by NSFC, no. 12122119 and no. 12341102. 

\section{Hyperbolic tetrahedra with degeneration}
\label{sec-31type}
\subsection{Type $3$-$1$ case}
In this section, we will review the concepts and results about type $3$-$1$ hyperbolic tetrahedra. More details about type $3$-$1$ hyperbolic tetrahedra could be found in \cite{Feng2023}.

The following is an description of $3$-$1$ type hyperbolic tetrahedra. We denote
$\IK^3 \subset \R^3$ by
the open ball representing $\IH^3$  via the Klein model. Then we can obtain a $3$-$1$ type
hyperbolic tetrahedron by the following process. Let $\mathscr{P}\subset \R^3$ be a convex
Euclidean tetrahedron
such that one vertex $v_1,$ lies on
the boundary of $\IK^3$  and the other vertices $v_i$, $i \in \{2,3,4\}$, lie in $\R^3\backslash \IK^3$, and each edge from $v_i$ intersects $\IK^3$. Let $C_i$ be the cone with
the apex $v_i$ tangent to $\partial\IK^3$ and $\pi_i$ be the half-space not containing $v_i$
such that $\partial \pi_i \cap \partial \IK^3= C_i \cap \partial \IK^3$. Then a $3$-$1$ type
hyperbolic tetrahedron is given by $\sigma= \mathscr{P}\cap (\cap_{i = 2}^4  \pi_i)$, and $v_1$ is called the
ideal vertex of $\sigma$, and the other vertices $v_2, v_3, v_4$ are called hyper-ideal vertices of $\sigma$.
$\Delta_i = \sigma\cap \partial \pi_i$ is called the vertex triangle of the $3$-$1$ type
hyperbolic tetrahedrons $\sigma$ corresponding to $v_i$, $2\leq i\leq 4$. By definition, a decorated $3$-$1$ type hyperbolic tetrahedron is a pair $(\sigma, \{H_1\})$, where $\sigma$ is a $3$-$1$ type hyperbolic tetrahedron defined above and $H_1$ is a horosphere centered at the ideal vertex $v_1$. Since $H_1$ is a horosphere, the intersection $\Delta_1=\sigma\cap H_1$ is an Euclidean triangle, which is also called a vertex triangle (corresponding to the ideal vertex $v_1$).
Note that, unlike $\Delta_1$, which is Euclidean, the three vertex triangles $\Delta_i$, $2\leq i\leq 4$ are all hyperbolic. Sometimes, to emphasize the difference between $\Delta_1$ and the other three $\Delta_i$, we also call $\Delta_1$ an ideal vertex triangle. We remark that in article \cite{FengGLiu}, we denote 1, 2, 3 by the three hyper-ideal vertices and 4 as the ideal vertex. This is not the same as the notation in this article. In this article, the vertex 1 is labeled as an ideal vertex and vertices 2, 3, 4 are labeled as hyper-ideal vertices.

Given a decorated $3$-$1$ type hyperbolic tetrahedron $(\sigma, \{H_1\})$, denote $e_{pq}$ by the side connecting $\Delta_p$ and $\Delta_q$ (in general, we don't distinguish between a vertex $i$ and the vertex triangle $\Delta_i$, so $e_{pq}$ can also be seen as the side connecting $p$ and $q$).
Let $\{p,q,r\}\subset\{1,2,3,4\}$, we denote $H_{pqr}$ by the polygonal face adjacent to $e_{pq}$, $e_{pr}$ and $e_{qr}$. The vertex edge $\Delta_p\cap H_{pqr}$ is denoted by $x^p_{qr}$ and the dihedral angle at $e_{pq} $ is denoted by $\alpha_{pq}$. Moreover, set $\alpha_{pq}= \alpha_{qp}$ for convenience. Note each dihedral angle between every hyperbolic polygonal face and the vertex triangle is always $\pi/2$. The moduli space of all dihedral angles of decorated $3$-$1$ type tetrahedron is
denoted by
\[
    \CB_{3|1} =\bigl\{ (\alpha_{12},\ldots,\alpha_{34}) \in {\R}^6_{>0} \mid \sum_{j\neq 1} \alpha_{1j} = \pi,  \sum_{j\neq i} \alpha_{ij} < \pi \; \text{for}\; i\in \{2,3,4\}   \bigr\} .
\]

For each pair $\{p,q\}\subset\{1,2,3,4\}$, denote $l_{pq}$ by the (signed) edge length of $e_{pq}$. Note $l_{23}$, $l_{34}$ and $l_{24}$ are always positive, while $l_{12}$, $l_{13}$ or $l_{14}$ maybe negative or zero, and $l_{1i}$ ($i\in\{1,2,3\}$) is negative (zero, resp.) exactly when $\Delta_i$ is inside (tangent to, resp.) the horosphere $H_1$. $l_{pq}$ is also called the decorated edge length, and
$$l = (l_{12},l_{13},l_{14},l_{23},l_{24},l_{34})\in {\R}^3 \times {\R}^3_{>0}$$
is called the decorated edge length vector. For the decorated $3$-$1$-type tetrahedron $(\sigma, H_1)$ with the decorated edge length vector $l$,
let $\theta_{ij}^k$ be the length of $x^k_{ij}$ which is the vertex edge or the intersection part of horospheres $H_1$ and the face $H_{1ij}$,  we define a function $\phi$ as follows:

\begin{lem}[Definition 3.4 and Lemma 3.4, \cite{Feng2023}]
\label{angleedge31}
For each pair $\{p,q\}\subset\{1,2,3,4\}$ and each vector $l\in {\R}^3 \times {\R}^3_{>0}$, define $\phi_{pq}(l)=\phi_{qp}(l)$, and
\begin{equation}\label{angle1}
  \phi_{1i}(l)=
  \frac{e^{l_{1j}}e^{l_{1k}}+e^{l_{1i}}(e^{l_{1k}}c_{ij} + e^{l_{1j}}c_{ik}-e^{l_{1i}}c_{jk})}{[(e^{2l_{1j}}+2e^{l_{1j}+l_{1i}}c_{ij}+e^{2l_{1i}})(e^{2l_{1k}}+2e^{l_{1k}+l_{1i}}c_{ik}+e^{2l_{1i}})]^{1/2}};
\end{equation}

\begin{equation}\label{angle2}
\phi_{ij}(l)=\frac{(e^{l_{1j}}+e^{l_{1i}} c_{ij})(c_{jk} + c_{ij} c_{ik})-(e^{l_{1k}}+e^{l_{1i}}c_{ik})s^2 _{ij}}{[(c^2 _{jk} + c^2_{ij} + c^2_{ik} + 2c_{jk}c_{ij}c_{ik}-1)(e^{2l_{1j}}+2e^{l_{1j}+l_{1i}}c_{ij}+e^{2l_{1i}})]^{1/2}};
\end{equation}
for $\{i,j,k\}=\{2,3,4\}$, where $c_{ij} = \cosh l_{ij}$, $s_{ij} = \sinh l_{ij}$. Then the right hand side of (\ref{angle2}) is symmetry about the index $i$ and $j$, which matches the setting $\phi_{pq}(l)=\phi_{qp}(l)$ at the beginning. In addition, if $l$ is the decorated edge length vector of $(\sigma, \{H_1\})$, then $ \phi_{ij}(l) = \cos \alpha_{ij}(l)$.
\end{lem}

\begin{prop}
\label{prop-extend-to-closer}
Each function $\phi_{pq}(l)$ ($1\leq p<q\leq4$) extends continuously to ${\R}^3 \times {\R}^3_{\geq0}$, which is the closer of ${\R}^3 \times {\R}^3_{>0}$, and the function form remains unchanged, still defined by equations (\ref{angle1}) and (\ref{angle2}). Furthermore, for $\{i,j\} \subset\{2,3,4\}$,
$\phi_{ij}(l) = 1$ when $ l_{ij} = 0 $.
\end{prop}

Proposition \ref{prop-extend-to-closer} may be derived easily from (\ref{angle1}) and (\ref{angle2}), hence we omit the details.

\subsection{Type $4$-$0$ case}
The geometry and degenerative behavior of type $4$-$0$ tetrahedra had been analyzed and studied in detail in \cite{Luo2018} and \cite{Feng2022}. A decorated $4$-$0$ type tetrahedron $\sigma$ is a (strictly) hyper-ideal hyperbolic tetrahedron (with no decorations since $\sigma$ has no ideal vertices). Similar conclusions like Lemma \ref{angleedge31} and Proposition \ref{prop-extend-to-closer} had been obtained in \cite{Luo2018} and \cite{Feng2022}. Here we only sketch some differences.
The moduli space of all dihedral angles of the $4$-$0$ type tetrahedron is
\[
    \CB_{4|0} =\bigl\{ (\alpha_{12},\ldots,\alpha_{34}) \in {\R}^6_{>0} \mid \sum_{j\neq i} \alpha_{ij} < \pi \; \text{for}\; i\in \{1, 2,3,4\}   \bigr\} .
\]
Since there is no ideal vertex, the decorated edge length is always positive, i.e.
$$l = (l_{12},l_{13},l_{14},l_{23},l_{24},l_{34})\in {\R}^6_{>0}.$$
For each pair $\{p,q\}\subset\{1,2,3,4\}$ and each vector $l\in {\R}^6_{>0}$, the $\phi$-functions $\phi_{pq}(l)$ are introduced (Lemma 4.3, \cite{Luo2018}), symmetric about the index $p$ and $q$, and satisfy $\phi_{pq}(l) = \cos \alpha_{pq}(l)$ for each decorated edge length vector $l$. Moreover, the function $\phi_{pq}$ extends continuously to $l\in {\R}^6_{\geq0}$, and $\phi_{pq}(l) = 1$ when $ l_{pq}=0$ (Lemma 4.7, \cite{Luo2018}).

\subsection{The cov-fucntion}
Now fix a decorated hyperbolic tetrahedron $\sigma$ in this subsection, either of type $3$-$1$, or of type $4$-$0$. Denote
$$\mathcal{R}={\R}^3 \times {\R}^3_{>0}$$
for type $3$-$1$ case and
$$\mathcal{R}={\R}^6_{>0}$$
for type $4$-$0$ case, respectively. It seems that by defining $\alpha_{ij}(l)=\cos^{-1}\phi_{ij}(l)$, the definition of dihedral angle $\alpha_{ij}$ can be extended continuously to any $l\in\overline{\mathcal{R}}$. However, the value of $\phi_{ij}(l)$ here may be greater than $1$ or less than $-1$. Therefore, $\alpha$ cannot be simply extended in this way.

\begin{prop}[Proposition 3.7, \cite{Feng2023}]
\label{realmetric}
The space of all hyperbolic tetrahedra $\sigma$ parameterized by the decorated edge length vector is
\[
  \CL = \{(l_{12},\ldots,l_{34}) \in \mathcal{R}\mid \phi_{ij}(l)\in(-1,1) \; \text{for all } \; \{i,j\} \subseteq \{1,2,3,4\}\},
\]
where $\mathcal{R}$ takes the from of ${\R}^3 \times {\R}^3_{>0}$ or ${\R}^6_{>0}$ depending on whether to deal with type $3$-$1$ or type $4$-$0$ tetrahedra.
\end{prop}

\begin{defi}
\label{defi:dihe}
Let $l\in \overline{\mathcal{R}}$. For any $\{i,j\}\subset \{1,2,3,4\}$, define
$$\alpha_{ij}(l)= \arccos (\min\{1,\max\{\phi_{ij}(l),-1\}\}).$$
$\alpha_{ij}$ is called the \textbf{extended dihedral angle} at the edge $e_{ij}$, and $l$ is called a generalized decorated edge length vector (or \textbf{generalized length} for short). A $3$-$1$ type ($4$-$0$ type, resp.) \textbf{generalized tetrahedron} $(\sigma,l)$ is a topological tetrahedron $\sigma$ which is combinatorially equivalent to some decorated $3$-$1$ type ($4$-$0$ type, resp.) tetrahedron, and a generalized length $l\in\overline{\mathcal{R}}$ which endows each edge $e_{ij}$ of $\sigma$ a number $l_{ij}$ as the decorated edge length. Sometimes for emphasis, we also refer to the previously appeared hyperbolic tetrahedron as a \textbf{real tetrahedron}, because it can be realized in $\IK^3 \subset \R^3$. In other words, $(\sigma,l)$ a real tetrahedron if and only if $l\in\CL$. In contrast, $(\sigma,l)$ is called a \textbf{degenerate tetrahedron} if $l\notin\CL$.
\end{defi}

We remark that there is no one or two dimensional degeneracy for each $l\in \mathcal{R}$. More specifically, consider the type $3$-$1$ case, the three pentagons corresponding to $v_1v_iv_j$ ($\{i, j\}\subset\{2, 3, 4\}$), and the triangle corresponding to $v_2v_3v_4$ (with corresponding decorations and truncations) always exist, and hence those $\{\theta^r_{pq}\}$ ($\{p, q, r\}\subset\{1, 2, 3, 4\}$) are always meaningful. See Theorem 2.2 in \cite{Feng2023} for more explanations. For the case of $4$-$0$ type, things are similar.

\begin{defi}
\label{def-l-l+}
For each $l=(l_{12}, l_{13}, l_{14}, l_{23}, l_{24}, l_{34}) \in {\R}^6$, set
$$l^+=(l_{12},l_{13}, l_{14}, l_{23}^+, l_{24}^+, l_{34}^+) \in {\R}^3\times {\R}^3_{\ge 0}=\overline{\mathcal{R}}$$
for type $3$-$1$ case, and set
$$l^+=(l_{12}^+, l_{13}^+, l_{14}^+, l_{23}^+, l_{24}^+, l_{34}^+) \in {\R}^6_{\ge 0}=\overline{\mathcal{R}}$$
for type $4$-$0$ case, where $ l_{ij}^+ = \max {\{0, l_{ij}\}}$.
\end{defi}

Then we could extend the definition of each dihedral angle $\alpha_{ij}$ continuously from the open set $\CL\subset\mathcal{R}$ to the whole space $\R^6$. First extend $\alpha_{ij}(l)$ from $l\in\CL$ to $l\in\overline{\mathcal{R}}$ by definition \ref{defi:dihe}, and then define $\alpha_{ij}(l)=\alpha_{ij}(l^+)$ for each $l\in\R^6$ according to Definition \ref{def-l-l+}. Both the extensions in the two steps are all continuously. Therefore we obtain a new continuous 1-form $\mu$ on ${\R}^6$ by
\[
  \mu(l) = \sum_{i\neq j} \alpha_{ij}(l) \di l_{ij},
\]
which is originally defined on $\CL$. Obviously $\mu(l)$ is smooth on $\CL$. Moreover, it is closed on $\CL$ due to the Schl\"{a}fli formula $2\text{d}vol+\sum_{i\neq j}l_{ij}d\alpha_{ij}=0$. Since $\CL$ is open and simply connected, one get a smooth co-volume function $cov:\CL\rightarrow \mathbb{R}$ by
$$cov_{\sigma}(l)=\int^l\mu.$$
We refer to \cite{Feng2023}, \cite{Feng2022} and \cite{Luo2018} for more details.
\begin{prop}[Proposition 4.7 in \cite{Feng2023} and Proposition 4.10 in \cite{Luo2018}]
\label{closed}
The continuous differential 1-form $\mu(l)=\sum\limits_{i\neq j} {\alpha}_{ij}(l) \,\mathrm{d} l_{ij}$ is closed in $\R^6$, i.e. for any Euclidean triangle $\Delta$ in $\R^6$, $\int_{\partial \Delta} \mu = 0$.
\end{prop}

\begin{cor}[Corollary 4.8 in \cite{Feng2023} and Corollary 4.12 in \cite{Luo2018}]
\label{cov}
The function $cov : \R^6 \to \R $ defined by the integral
  \[cov_{\sigma}(l)= \int_{(0,\ldots,0)}^{l} \mu + cov(0,\ldots,0) \]
is a $C^1$-smooth convex function.
\end{cor}

\section{Some geometric properties of type $3$-$1$ and $4$-$0$ tetrahedra}
\label{sec-proper}
\subsection{Some geometric properties of type $3$-$1$ tetrahedra}
\label{sec-31proper}
Given a general metric $l\in\R^3 \times \R^3_{\ge 0}$, recall Proposition \ref{prop-extend-to-closer} shows that $\phi_{ij}(l)$ is well defined by the formula (\ref{angle1}) and (\ref{angle2}). In the following, we set
$$x_2=e^{l_{12}}, \;x_4=e^{l_{13}}, \;x_6=e^{l_{14}}, \;x_1=\cosh{l_{23}}, \;x_3=\cosh{l_{34}}, \;x_5=\cosh{l_{24}}.$$
Since the decorated length of hyper-ideal edges should be bigger than $0$, and the decorated length of ideal edges belongs to $\R$, so $x_1$, $x_3$, $x_5\ge1$ and
$x_2$, $x_4$, $x_6>0$. Particularly, any one of $x_1$, $x_3$, $x_5$ equals to $1$ means that the length of the corresponding hyper-ideal edge is $0$.

According to the law of cosine, if 
$\Delta_1$ is assumed to be an equilateral triangle with side length $1$, then in the decorated $3$-$1$ type hyperbolic tetrahedron, $x= (x_1, \dots, x_6)$ satisfies the following relation:
\begin{equation}\label{eq}\left\{\begin{array}{l}
x_1 = 2x_2x_4-\dfrac{1}{2}(\dfrac{x_2}{x_4}+\dfrac{x_4}{x_2}),\\[10pt]
x_3 = 2x_4x_6-\dfrac{1}{2}(\dfrac{x_4}{x_6}+\dfrac{x_6}{x_4}),\\[10pt]
x_5 = 2x_2x_6-\dfrac{1}{2}(\dfrac{x_2}{x_6}+\dfrac{x_6}{x_2}).
\end{array}\right.\end{equation}

Now, we introduce the following useful lemma:
\begin{lem}\label{bdd}
If $(x_1,x_2,x_3,x_4,x_5,x_6) \in \R_{\geq 0}^6$ satisfies (\ref{eq}), and $x_1$, $x_3$, $x_5 \geq 1$, then
\begin{enumerate}[(i)]
\item $x_2$, $x_4$, $x_6>\dfrac{1}{2},$
\item if $x_1$, $x_3$, $x_5 \leq A $, then $\dfrac{2(A+1)+\sqrt{2(A+1)}}{4A+2} \leq x_2$, $x_4$, $x_6 \leq \dfrac{A + \sqrt {A^2+3}}{3} \leq A.$
\end{enumerate}
\end{lem}
\begin{proof}

The three equations in (\ref{eq}) have the same form:
\[k = 2mn-\frac{1}{2}(\frac{m}{n}+\frac{n}{m}).\]
For $(i)$:

Since $x_1, x_3,x_5 \geq 1$,
\[1 \leq 2mn-\frac{1}{2}(\frac{m}{n}+\frac{n}{m}),\]
i.e. $$0 < 4m^2n^2-m^2-n^2-2mn .$$

Hence, $4m^2-1 > 0$ and $4n^2-1 > 0$, so $m,n > 1/2$. This means $x_2, x_4, x_6 > 1/2 $.

For $(ii)$:
If $x_1,x_3,x_5 \in [1, A], $ then
\[1 \leq 2mn-\frac{1}{2}(\frac{m}{n}+\frac{n}{m}) \leq A,\]
i.e. $0 < 4m^2n^2-m^2-n^2-2mn $ and  $4m^2n^2-m^2-n^2-2Amn < 0.$
Let
$$f(m,n) = 4m^2n^2 - m^2 - n^2 - 2mn. $$

Firstly we claim that, if $f(m,n) > 0$, then $m, n$ can not less than $1$ at the same time.

If $m, n < 1$, then $$f(m,n) = m^2n^2 - m^2 +m^2n^2 - n^2 + 2m^2n^2 - 2mn < 0.$$
This is in contradiction with  $0 < 4m^2n^2-m^2-n^2-2mn $, so the claim is hold.

Next, let
\[g(m,n) = 2mn-\frac{1}{2}\Big(\frac{m}{n}+\frac{n}{m}\Big).\]

Since $m,n >1/2$,
\[\frac{\partial g}{\partial m} =2n-\frac{1}{2n}+\frac{n}{2m^2} > 0,~~ \frac{\partial g}{\partial n} =2m-\frac{1}{2m}+\frac{m}{2n^2} > 0.\]

This means $g(m,n)$ increases monotonically with $m,n$. Without loss of generality, we can assume $ x_2 \geq x_4 \geq x_6,$ at this time, we have $ x_1 \geq x_5 \geq x_3,$ and we also have $x_2 \geq x_4 \geq 1$ by the claim above.

We can consider $4m^2n^2-m^2-n^2-2kmn = 0$ as an equation for $m$, then
\[m = n \cdot \frac{k + \sqrt{k^2+4n^2-1}}{4n^2 - 1}.\]

If $x_1 = A, $ since $x_2 > 1/2$, $4x_2^2 - 1 > 0$, then
\[x_4 = x_2 \cdot \frac{A + \sqrt{A^2+4x_2^2-1}}{4x_2^2 - 1}.\]

Based on our assumption, we have know that $x_4 \ge 1$, thus
\[x_2\Big(A + \sqrt{A^2+4x_2^2-1}\Big) \geq 4x_2^2-1. \]
\[(4x_2^2-1)(3x_2^2-2Ax_2-1) \le 0.\]
So
\[x_2 \leq \frac{A + \sqrt {A^2+3}}{3} \leq A.\]

On the other hand, if $x_3 = 1$, then \[x_4 =  \frac{x_6}{2x_6 - 1}.\] Since $g(m,n)$ increases monotonically with $m,n$,  we have \[g(x_4, x_4) = \frac{2x_6^2}{(2x_6 - 1)^2} - 1 \le g(x_2, x_4) = x_1 \le A.\] Thus, \[x_6 \ge \frac{2(A+1)+\sqrt{2(A+1)}}{4A+2}.\]

\end{proof}

Consider the map $G:(x_2,x_4,x_6) \to (x_1,x_3,x_5)$ defined by the equation (\ref{eq}), then we have the following result:

\begin{prop}\label{inv}
Let $\Omega^*$ be the region $(1/2, +\infty)^3 $ and the Jacobian matrix of $G$ is denoted as $J(G)$. Then, the determinant of $J(G)$ is positive on $\Omega^*$. Further, $G$ is invertible on $\Omega^*$.
\end{prop}
\begin{proof}
Firstly, we claim that $G$ is injective in $\Omega^*$.

If we assume there are two triples $(x_2,x_4,x_6)$ and $(\hat x_2, \hat x_4, \hat x_6)$ satisfying
$$G(x_2,x_4,x_6) = G(\hat x_2, \hat x_4, \hat x_6)= (x_1, x_3, x_5).$$

Without loss of generality, we can assume $x_2 \geq \hat x_2$. Since $g(m,n) = 2mn-\frac{1}{2}(\frac{m}{n}+\frac{n}{m})$ is increases monotonically with $m,n$, we have that $\hat x_4 \geq x_4$ and $\hat x_6 \geq x_6$. On the other hand, we also have $\hat x_6 \leq x_6$ from $\hat x_4 \geq x_4$. Thus, $\hat x_6 = x_6$. This implies that $\hat x_4 = x_4$ and $\hat x_2 = x_2$; that is say our claim is true.

Secondly, by directly calculation, $J(G)$  is
\begin{equation}\left[ \begin{array}{ccc}
2x_4-\dfrac{1}{2x_4}+\dfrac{x_4}{2x_2^2} & 2x_2-\dfrac{1}{2x_2}+\dfrac{x_2}{2x_4^2} & 0\\
0&                                       2x_6-\dfrac{1}{2x_6}+\dfrac{x_6}{2x_4^2} & 2x_4-\dfrac{1}{2x_4}+\dfrac{x_4}{2x_6^2} \\
2x_6-\dfrac{1}{2x_6}+\dfrac{x_6}{2x_2^2} & 0                                      & 2x_2-\dfrac{1}{2x_2}+\dfrac{x_2}{2x_6^2}
\end{array}\right]
\end{equation}
Since the determinant for matrix
$\left[ \begin{array}{ccc}
a & b & 0\\
0 & c & d \\
e & 0 & f
\end{array}\right]$
equals $acf+bde$, so this determinant must be positive if every element in this matrix is positive. While, we have shown that every element of $J(G)$ is positive previously. Thus $\det (J(G)) > 0$ on $\Omega^*$ and $G$ is invertible on $\Omega^*$.
\end{proof}

\begin{prop}
\label{surjection}
Given any $A \ge 1$, For any $(x_1, x_3,x_5) \in [1, A]^3$, there exists one and only one solution $(x_2, x_4,x_6) \in (1/2, A]^3$ to the equation (\ref{eq}). Thus, for $(x_1, x_3,x_5) \in [1, +\infty)^3$, we know the equation (\ref{eq}) of $(x_2, x_4,x_6)$ has a unique solution.
\end{prop}
\begin{proof}
For any given $x_1, x_3, x_5 \in [1,A]$, by Lemma \ref{bdd} and Proposition \ref{inv}, there is at most one triple $(x_2,x_4,x_6)$ satisfying the formula $(\ref{eq})$. Since
$A\geq 1$ is arbitrary chosen, we know that for any $(x_1, x_3,x_5)$, there is at most one solution of the equations (\ref{eq}), thus we just need to show the existence of the solution. Let $x_1 = X$, $x_3 =Y$, $x_5 = Z$, w.l.o.g., suppose $1\le X, Y, Z\le A$. Let's $ x_2 = a $ to be determined, then
\[x_4 = x_2 \cdot \frac{X + \sqrt{X^2 + 4x_2^2 - 1}}{4x_2^2 - 1} = a \cdot \frac{X + \sqrt{X^2+4a^2-1}}{4a^2 - 1},\]
\[x_6 = x_2 \cdot \frac{Z + \sqrt{Z^2 + 4x_2^2 - 1}}{4x_2^2 - 1} = a \cdot \frac{Z + \sqrt{Z^2+4a^2-1}}{4a^2 - 1}.\]
Let
\[g(m,n) = 2mn-\frac{1}{2}\Big(\frac{m}{n}+\frac{n}{m}\Big),\]
then
\begin{eqnarray*}
g(a)&\triangleq& g(x_4(a), x_6(a)) \nonumber\\
    &=& \frac{1}{2}\bigg[4a^2 \cdot \frac{(X + \sqrt{X^2+4a^2-1})(Z+ \sqrt{Z^2+4a^2-1})}{(4a^2 - 1)^2}\\
    &&  \quad -\frac{Z + \sqrt{Z^2+4a^2-1}}{X + \sqrt{X^2+4a^2-1}}- \frac{X+ \sqrt{X^2+4a^2-1}}{Z+ \sqrt{Z^2+4a^2-1}}\bigg]. \nonumber
\end{eqnarray*}
Thus, to prove this proposition, we just need to show that the range of $g(a)$ covers $[1, A]$ if $a$ runs over the interval $(\frac{1}{2}, A]$.
Set $k = 4a^2-1$, since $a \in (\frac{1}{2}, A]$, then $k \in (0, 4A^2-1]$, and
\begin{eqnarray*}
h(k)&\triangleq& g(a) \nonumber\\
    &=& \frac{1}{2}\bigg[(k+1) \cdot \frac{(X+ \sqrt{X^2+k})(Z+ \sqrt{Z^2+k})}{k^2}-\frac{Z+ \sqrt{Z^2+k}}{X+ \sqrt{X^2+k}}-\frac{X+ \sqrt{X^2+k}}{Z+ \sqrt{Z^2+k}}\bigg]. \nonumber
\end{eqnarray*}
Since $X, Z\in [1, A]$, we know that $g(a) = h(k) \to +\infty$ as $a \to 1/2$ ~(i.e. $k \to 0$). On the other hand, if $a=A$ (i.e. $k = 4A^2-1$), we have

\begin{eqnarray*}
g(A)= h(4A^2-1) &=& \frac{1}{2}\bigg[\frac{4A^2(X+\sqrt{X^2+4A^2-1})(Z+\sqrt{Z^2+4A^2-1})}{(4A^2-1)^2}\\
    &&  \quad -\frac{Z+ \sqrt{Z^2+4A^2-1}}{X+ \sqrt{X^2+4A^2-1}}-\frac{X+ \sqrt{X^2+4A^2-1}}{Z+ \sqrt{Z^2+4A^2-1}}\bigg] \nonumber\\
    &\le& \frac{1}{2}\bigg(\frac{4A^2 (A+ \sqrt{5A^2-1})^2}{(4A^2-1)^2} - 2\bigg).
\end{eqnarray*}
Since
\begin{eqnarray*}
&&\frac{4A^2 (A+ \sqrt{5A^2-1})^2}{(4A^2-1)^2} =\frac{4\big(1+\sqrt{1+(4-A^{-2})}\big)^2}{(4-A^{-2})^2}\\
&=&\frac{4(1+\sqrt{1+x})^2}{x^2}= 4[x^{-1}+\sqrt{x^{-2}+x^{-1}}]^2,
\end{eqnarray*}
where $x= 4-A^{-2}$. It's monotonically decreasing with $A$. Thus, for any $A\ge 1$,
$$g(A) \le 1.$$

By Intermediate Value Theorem, we know that the range of $g(a)$ covers $[1,A]$ for $a \in (\frac{1}{2}, A]$. Thus, for $x_3=Y\in [1, A]$, there is a solution of $g(a)=Y$. Further, for any $(x_1, x_3, x_5) \in [1, A]^3$, there is a solution $(x_2, x_4,x_6) $ of the equations (\ref{eq}).
\end{proof}

\begin{cor}
\label{prop-equa-length-1}
Every decorated $3$-$1$ type hyperbolic tetrahedron is uniquely determined (up to isometry) by the lengthes of its three hyper-ideal edges, when assuming that its ideal vertex triangular is a equilateral triangle with length $1$.
\end{cor}

In Corollary \ref{prop-equa-length-1}, the dihedral angles of the edges relative to the ideal vertex are all $\pi/3$. Thus, we only need to consider the change of the remaining three dihedral angles. This means that we only need to discuss the range of $\phi$ on the hyper-ideal edges. Further, due to rotation symmetry, we just need to consider the case where the hyper-ideal edge is $e_{23}$. By (\ref{angle2}), we have
\begin{eqnarray}
\phi_{23}(l)&=&\frac{(e^{l_{13}}+e^{l_{12}}c_{23})(c_{34}+c_{23}c_{24})-(e^{l_{14}}+e^{l_{12}}c_{24})s^2_{23}}{\sqrt{c^2_{23}+c^2_{24}+c^2_{34}+2c_{23}c_{24}c_{34}-1}\sqrt{e^{2l_{12}}+e^{2l_{13}}+2e^{l_{12}+l_{13}}c_{23}}} \notag \\
&=&\frac{(x_4+x_1x_2)(x_3+x_1x_5)-(x_6+x_2x_5)(x_1^2-1)}{[(x_1^2+x_3^2+x_5^2+2x_1x_3x_5-1)(x_2^2+x_4^2+2x_1x_2x_4)]^{1/2}} \notag\\
&=&\frac{x_4x_3+x_2x_5+x_1x_2x_3+x_1x_4x_5-x_1^2x_6+x_6}{[(x_1^2+x_3^2+x_5^2+2x_1x_3x_5-1)]^{1/2}[(x_2^2+x_4^2+2x_1x_2x_4)]^{1/2}}
\end{eqnarray}

Now we define $\phi (x_1, \dots, x_6) \triangleq \phi_{23}(l)$, which is a function of $(x_1, \dots, x_6)$, the following two results are obvious:
\begin{equation}\label{sym}
\phi(x_1,x_2,x_3, x_4,x_5,x_6)=\phi(x_1,x_4,x_5, x_2,x_3,x_6),
\end{equation}
and
\begin{equation}
\frac{\partial \phi}{\partial x_6} < 0.
\end{equation}
Let
\[A_0= \frac{1}{[(x_1^2+x_3^2+x_5^2+2x_1x_3x_5-1)]^{1/2}[(x_2^2+x_4^2+2x_1x_2x_4)]^{3/2}},\]
\[A_1= \frac{1}{[(x_1^2+x_3^2+x_5^2+2x_1x_3x_5-1)]^{3/2}[(x_2^2+x_4^2+2x_1x_2x_4)]^{1/2}},\]

then
\begin{equation}\label{eq:pde}\left\{\begin{array}{l}
\dfrac{\partial \phi}{\partial x_2} = A_0(x_1^2-1)[x_4(x_1x_6 + x_2x_3 - x_4x_5) + x_2x_6], \\[9pt]
\dfrac{\partial \phi}{\partial x_2} = A_1(x_1^2-1)[x_5(x_1x_6 + x_2x_3 - x_4x_5) + x_4 + x_1x_2 +x_3x_6],\\[9pt]
\dfrac{\partial \phi}{\partial x_4} = A_0(x_1^2-1)[x_2(x_1x_6 - x_2x_3 + x_4x_5) + x_4x_6],\\[9pt]
\dfrac{\partial \phi}{\partial x_5} = A_1(x_1^2-1)[x_3(x_1x_6 - x_2x_3 + x_4x_5) + x_2 + x_1x_4 +x_5x_6],
\end{array}\right.\end{equation}

\begin{prop}\label{prop:twop0}  For any $x = (x_1, \dots, x_6) \in \R^6_{> 0}$, then at least one of the following holds:
\begin{enumerate}[(i)]
\item $\dfrac{\partial \phi}{\partial x_2}>0$ and $\dfrac{\partial \phi}{\partial x_3}>0$;
\item $\dfrac{\partial \phi}{\partial x_4}>0$ and $\dfrac{\partial \phi}{\partial x_5}>0$.
\end{enumerate}
\end{prop}

\begin{proof}
By \eqref{eq:pde}, if $x_2x_3 - x_4x_5 < 0$, then $\frac{\partial \phi}{\partial x_4}>0, \frac{\partial \phi}{\partial x_5}>0$, so the assertion $(ii)$ holds. Otherwise, the assertion $(i)$ holds.
\end{proof}

Next, we use the partial derivative $\frac{\partial \phi}{\partial x_i}$ to perform the monotonicity analysis, $i= 2,3,4,5.$ For any such $i,$ we define the index $\hat{i},$ which is in pair of $i,$ as follows,

\begin{equation*}\hat{i}=\left\{\begin{array}{cc}
3, & i=2,\\
2, & i=3,\\
5, & i=4,\\
4, & i=5.
\end{array}\right.\end{equation*}

For any $k\in \{1,\cdots, 6\},$
we denote $m_k=(0,\cdots, 1,\cdots, 0)$ by  the $k$-th coordinate unit vector in $\R^6.$

\begin{prop}\label{prop:twop}
Let $x\in \R^6_{> 0}$ , $x_1 > 1$, $i\in \{2,3,4,5\}.$
\begin{enumerate}[(i)]
\item If $\dfrac{\partial \phi}{\partial x_i}(x)\geq 0,$ then for any $t,s\geq 0,$ $$\frac{\partial \phi}{\partial x_i}(x+tm_i+sm_{\hat{i}})\geq 0.$$ In particular,
$\phi(x)\leq \phi(y),$ for any $y$ satisfying $y_i\geq x_i, y_j=x_j$ ($\forall j\neq i$).
\item If $\dfrac{\partial \phi}{\partial x_i}(x)\leq 0,$ then for any $t,s\leq 0$ satisfying $x+tm_i+sm_{\hat{i}}\in \R^6_{> 0},$  $$\frac{\partial \phi}{\partial x_i}(x+tm_i+sm_{\hat{i}})\leq 0.$$
In particular,
{$\phi(x)\leq \phi(y),$} for any $y$ satisfying $ y_i\leq x_i, y_j=x_j$  ($\forall j\neq i$).
\end{enumerate}
\end{prop}

\begin{proof}
Without loss of generality, we prove the result for $i = 2$ and assume that $\frac{\partial \phi}{\partial x_2}(x)\geq 0$.
By equations (\ref{eq:pde}) , since $A_0 (x_1^2-1)> 0$, the sign of $\frac{\partial \phi}{\partial x_2}(x)$ is decided by
\[G(x) = x_4(x_1x_6 + x_2x_3 - x_4x_5) + x_2x_6.\]
Obviously, if fix $x_1,x_4,x_5,x_6$, $G(x)$ is increasing monotonously with $x_2, x_3$. So, when $  G(x) \geq 0$, for any $t, s \geq 0$, the following is hold

\[G(x + t m_2 + s m_3 ) \geq G(x) \geq 0.\]

\end{proof}

Following this proposition, we can get a useful estimate.

\tm\label{thm:kest1}
Let $x\in \R^6_{> 0}$ satisfying $1 \leq x_3, x_5 \leq b, 1/2 \leq x_2,x_4, x_6 \leq a. $ Then
\begin{eqnarray*}
\phi(x) &\leq &\max \{\phi(x_1,a,b,a,b,1/2), \;\phi(x_1,a,b,a,1,1/2), \\
&& \quad\quad\;\phi(x_1,a,b,1/2,b,1/2), \; \phi(x_1,a,b,1/2,1,1/2)\}.
\end{eqnarray*}
\tmd

\begin{proof}
Since $\frac{\partial \phi}{\partial x_6} < 0,$

\[\phi(x)\leq \phi(x_1, x_2,x_3,x_4,x_5,1/2).\]

On the other hand, by Proposition \ref{prop:twop0} and (\ref{sym}), for $(x_1, x_2,x_3,x_4,x_5,1/2)$, we could assume that $\frac{\partial \phi}{\partial x_2}>0, \frac{\partial \phi}{\partial x_3}>0$, without loss of generality.

By Proposition \ref{prop:twop},
\[\phi(x)\leq \phi(x_1, x_2,x_3,x_4,x_5,1/2) \leq \phi(x_1,a,b,x_4,x_5,1/2).\]

\begin{description}
  \item[Case 1]
  If $\frac{\partial \phi}{\partial x_4} \geq 0$ at the point$(x_1,a,b,x_4,x_5,1/2)$, then
  \[\phi(x_1,a,b,x_4,x_5,1/2) \leq \phi(x_1,a,b,a,x_5,1/2). \]
  Further, consider  the following two cases at $(x_1,a,b,a,x_5,1/2)$.
\begin{description}
  \item[Case 1.1]
  If$\frac{\partial \phi}{\partial x_5} \geq 0$, then $\phi(x_1,a,b,a,x_5,1/2) \leq \phi(x_1,a,b,a,b,1/2);$
\item[Case 1.2]
  If $\frac{\partial \phi}{\partial x_5} < 0$, then $\phi(x_1,a,b,a,x_5,1/2) \leq \phi(x_1,a,b,a,1,1/2).$
\end{description}
  Thus \[\phi(x)\leq \max \{\phi(x_1,a,b,a,b,1/2), \phi(x_1,a,b,a,1,1/2)\};\]

 \item[Case 2]
 If $\frac{\partial \phi}{\partial x_4} < 0$ at the point $(x_1,a,b,x_4,x_5,1/2)$, then
\[\phi(x_1,a,b,x_4,x_5,1/2) \leq \phi(x_1,a,b,1/2,x_5,1/2). \]
Further, consider  the following two cases at $(x_1,a,b,1/2,x_5,1/2)$.
 \begin{description}
  \item[Case 2.1] If $\frac{\partial \phi}{\partial x_5} \geq 0$, then $\phi(x_1,a,b,1/2,x_5,1/2) \leq \phi(x_1,a,b,1/2,b,1/2);$
  \item[Case 2.2] If $\frac{\partial \phi}{\partial x_5} < 0$, then $\phi(x_1,a,b,1/2,x_5,1/2) \leq \phi(x_1,a,b,1/2,1,1/2).$

  Thus \[\phi(x)\leq \max \{\phi(x_1,a,b,1/2,b,1/2),\; \phi(x_1,a,b,1/2,1,1/2)\}.\]
\end{description}

\end{description}
In conclusion,
\begin{eqnarray*}
\phi(x)&\leq & \max \{\phi(x_1,a,b,a,b,1/2),\; \phi(x_1,a,b,a,1,1/2), \\
            &&\quad\quad\;\phi(x_1,a,b,1/2,b,1/2), \;\phi(x_1,a,b,1/2,1,1/2)\}.
\end{eqnarray*}
\end{proof}

By Theorem \ref{thm:kest1}, if $x_1 = k$ is the maximum in $\{x_1,x_3,x_5\}$, then
\begin{eqnarray*}
\phi(x)&\leq & \max \{\phi(k,k,k,k,k,1/2), \;\phi(k,k,k,k,1,1/2), \\
&&\quad\quad\; \phi(k,k,k,1/2,k,1/2),\; \phi(k,k,k,1/2,1,1/2)\}\\
       &\triangleq & \max \{h_1,h_2,h_3,h_4\}.
\end{eqnarray*}

\begin{equation}\label{eq:h}\left\{\begin{array}{l}
h_1 = \dfrac{4k^3+3k^2+1}{2\sqrt{(2k^3+2k^2-1)(2k^3+2k^2)}} ,\\[10pt]
h_2 = \dfrac{2k^3+3k^2+2k+1}{4k\sqrt{2k^3+2k^2}},\\[10pt]
h_3 = \dfrac{2k^3+2k^2+k+1}{\sqrt{(2k^3+3k^2-1)(8k^2+1)}},\\[10pt]
h_4 = \dfrac{2k^3-k^2+4k+1}{2k\sqrt{8k^2+1}}.
\end{array}\right.\end{equation}

In the following of this section, we abbreviate the lengthes $l_{23}$, $l_{12}$, $l_{34}$, $l_{13}$, $l_{24}$ and $l_{14}$ as $l_1$, $l_2$, $l_3$, $l_4$, $l_5$ and $l_6$ respectively. We relabel these $\{\phi_{ij}\}$ as $\phi_i$ ($1\leq i\leq 6$) in the same way. However, in Section \ref{sec-erf} and beyond, the notation $l_i$, with the lower index $i$ refers specifically to a certain marked edge with $i\in E=\{1,\cdots,m\}$, have different meanings.
\co\label{lem:longest}
Let $\sigma$ is a general $3$-$1$ hyperbolic tetrahedra with the ideal vertex $v_1$, if its decoration $H$ satisfies the ideal vertex triangle corresponding to $v_1$ is the length $1$ equilateral triangle, and $l_1=\max {\{l_1, l_3, l_5\}},$ then
\begin{enumerate}[(i)]
\item let $x_1=\cosh l_1, $\begin{eqnarray}\label{cor-bd}\cos\alpha_1 & \leq & \max\big\{\phi(x_1,x_1,x_1,x_1,x_1,1/2), \; \phi(x_1,x_1,x_1,x_1,1,1/2),\\&&\quad\quad\;\;\phi(x_1,x_1,x_1,1/2,x_1,1/2), \; \phi(x_1,x_1,x_1,1/2,1,1/2) \big\},\nonumber \end{eqnarray}

\item let $l_1=\arccosh 2,$  $$\alpha_1> \frac{2\pi}{11}.$$
\end{enumerate}
\cod
\begin{proof}
 By Lemma \ref{bdd}, if $l_1=\max {\{l_1, l_3, l_5\}},$ then $x_2,x_4,x_6 \le x_1$. Let $a = b = x_1$, we get (\ref{cor-bd}) from Theorem \ref{thm:kest1}.

  Further, if $l_1 =\max {\{l_1, l_3, l_5\}} =\arccosh 2$, i.e. $x_1 =\max {\{x_1, x_3, x_5\}} = 2$, by Lemma \ref{bdd} (2), we can get a better bound of $\frac{4}{5} < x_2,x_4,x_6 < \frac{8}{5}$. Let $a = \frac{8}{5}$ and $x_1 = b = 2$, and we replace the lower bound $\frac{1}{2}$ with $\frac{4}{5}$, then we can get a more accurate estimate by using the same proof as Theorem \ref{thm:kest1}.
\begin{eqnarray*}
\phi(2)&\leq & \max \{\phi(2,8/5,2,8/5,2,4/5), \;\phi(2,8/5,2,8/5,1,4/5), \\
    &&\quad\quad\; \phi(2,8/5,2,4/5,2,4/5),\; \phi(2,8/5,2,4/5,1,4/5)\}\\
       & = & \frac{7\sqrt{2}}{12}.
\end{eqnarray*}
 Thus, \[\alpha_1 \geq \arccosh{\frac{7\sqrt{2}}{12}} >\frac{2\pi}{11}.\]
\end{proof}

Now, we have the following theorem:
\begin{theo}\label{Range}
Let $\sigma$ be a general $3$-$1$ hyperbolic tetrahedra with the ideal vertex $v_1$, if its decoration $H_1$ satisfies that the vertex triangle corresponding to $v_1$ is a equilateral triangle with length $1$, if $0 < l_1$, $l_3$, $l_5 \leq \arccosh 2$, then
$$\phi_i(l)\in (-1,1),\quad  1\leq i\leq 6,$$
where $\{\phi_i\}$ are the relabeling of $\{\phi_{ij}\}$. Particularly, $l\in \mathcal{L}$, i.e. $\sigma$ is a real hyperbolic tetrahedron.
\end{theo}

\begin{proof}
We just need to prove $\phi_{23}(l), \phi_{24}(l), \phi_{34}(l) \in (-1, 1)$. By the rotation symmetry, we just prove $\phi_{23}(l) \in (-1, 1)$ here.

Since $0 < l_1$, $l_3$, $l_5 \leq \arccosh 2,$ then $x_1$, $x_3$, $x_5 \in (1,2]$ and $x_2$, $x_4$, $x_6 \in (1/2, 2]$ by Lemma \ref{bdd}.
By Theorem \ref{thm:kest1}
\[
\phi_{23}(x) \leq \max\{\tau_1(x_1), \tau_2(x_1), \tau_3(x_1), \tau_4(x_1)\},
\]
where
\[\tau_1(x)=\phi(x,2,2,2,2,1/2) = \frac{-x^2 + 16x + 17}{4\sqrt{(2+2x)(x^2+8x+7)}},\]
\[\tau_2(x)=\phi(x,2,2,2,1,1/2) = \frac{-x^2 + 12x + 13}{4(x+2)\sqrt{2+2x}},\]
\[\tau_3(x)=\phi(x,2,2,1/2,2,1/2) = \frac{-x^2 + 10x + 11}{\sqrt{(8x+17)(x^2+8x+7)}},\]
\[\tau_4(x)=\phi(x,2,2,1/2,1,1/2) = \frac{-x^2 + 9x + 7}{(x+2)\sqrt{8x+17}}.\]

Obviously, by direct calculation, the above four functions are monotonously decreasing when $x \geq 1$, and the value of every function is $1$ at $x=1$. Thus $\phi_{23}(l) < 1$.

On the other hand, since $\phi_{23}(x)$ is monotonously decreasing for $x_6$ , we can assume $x_6 = 2$. At the same time , we can suppose that $\phi_{23}(x) < 0$, otherwise the conclusion is valid. In these assumptions, the numerator of $\phi_{23}(x)$ should be negative, i.e.
\[x_2x_3+x_4x_5 < x_1x_6 \leq 2x_1.\]

Based on the equation (\ref{eq:pde}), $\frac{\partial \phi_{23}(x)}{\partial x_i} > 0$ ($i = 2,3,4,5$), this means $\phi_{23}(x)$ is monotonically increasing about $x_2,x_3,x_4,x_5$. So, if $\phi_{23}(x) < 0$ , for any $x_1 \in [1,2]$

\[\phi_{23}(x) \geq \phi(x_1,1/2,1,1/2,1,2) = -\frac{\sqrt{2}(2x_1-3)}{\sqrt{x_1+1}} \geq -\sqrt{\frac{2}{3}} > -1 \]

In conclusion, if $0 < l_1,l_3,l_5 \leq \arccosh 2,$ then $$\phi_i(l)\in (-1,1)\quad  1\leq i\leq 6.$$
\end{proof}

By the definition of $\phi_{ij}$ and $\alpha_{ij},$ they are continuous functions. And if fix $x_2,\cdots,x_6,$ $\phi$ goes to $1$ as $x_1$ goes to $0$, this means $\alpha_1$ goes to $0$. In the following, we give a quantitative estimate.
\begin{prop}\label{prop:c0}
For any $C > 1$ and $\epsilon > 0,$ there exists $\delta(C,\epsilon) > 0$ such that
for $$\max_{2\leq i\leq 6}\{x_i\}\leq C\quad \mathrm{and}\quad \ x_1\leq 1+\delta,$$ we have $$\phi(x)\geq \cos\epsilon.$$ In particular, $\alpha_1\leq \epsilon.$  For the case that $C = 2$ and {$\epsilon \ll \pi,$} we can choose $$\delta(2,\epsilon)=\frac{1}{2\pi^2}\epsilon^2.$$ \end{prop}

\begin{proof}
By Lemma \ref{bdd}, noting that $x_3,x_5\in [1,C], x_2,x_4,x_6 \in [1/2, C]$ and $x_1\in [1,1+\delta],$
\begin{equation}
\phi_{23}(x)\geq\frac{x_4x_3+x_2x_5+x_2x_3+x_4x_5-(2\delta+\delta^2)C}{[(x_1^2+x_3^2+x_5^2+2x_1x_3x_5-1)]^{1/2}[(x_2^2+x_4^2+2x_1x_2x_4)]^{1/2}}.
\end{equation}

Choose $\delta<\frac{1}{2C}$ such that the numerator is positive. Hence
\begin{eqnarray}
\phi(x) &\geq&
\frac{(x_2+x_4)(x_3+x_5) - (2\delta+\delta^2)C}{\sqrt {(x_2+x_4)^2 + 2\delta C^2}\sqrt {(x_3+x_5)^2 + 2\delta C^2 + 2\delta+\delta^2}}\nonumber\\
&=&\frac{1 - \frac{(2\delta+\delta^2)C}{(x_2+x_4)(x_3+x_5)}}{\sqrt {1 + (2\delta C^2)/(x_2+x_4)^2}\sqrt {1 + (2\delta C^2 + 2\delta+\delta^2)/(x_3+x_5)^2}}\nonumber\\
&\geq &\frac{1 - \frac{1}{2}(2\delta+\delta^2)C}{1 + (2\delta C^2 + 2\delta+\delta^2)}\to 1 \quad \quad (\delta\to 0).\label{eq:del1}
\end{eqnarray}
Thus, for any $\epsilon>0,$ there exists $\delta(C,\epsilon)>0$ s.t $\phi(x)\geq \cos\epsilon.$ From the definition of generalized dihedral angle, $\alpha_1\leq \epsilon.$

For the second assertion, for $C=2$ and {$\epsilon < \pi,$} we choose
$\delta=\frac{1}{6\pi^2}\epsilon^2<\frac{1}{2}.$

Then
\begin{eqnarray*}
\phi(x) &\geq&\frac{1 - (2\delta+\delta^2)}{1 + (10\delta+\delta^2)}\\
&\geq&\frac{1 - \frac32\delta}{1 + \frac{21}{2}\delta}\geq 1-\frac{2}{\pi^2}\epsilon^2.
\end{eqnarray*} Note that for any {$0 < \epsilon < \pi,$}
\begin{equation*}1-\cos\epsilon=2\sin^2(\frac{\epsilon}{2}) \geq \frac{2}{\pi^2}\epsilon^2.
\end{equation*}
Thus $\phi(x)\geq \cos\epsilon.$
\end{proof}

\subsection{Some geometric properties of type $4$-$0$ tetrahedra}
\label{sec-40proper}
By the results proved earlier for a type $4$-$0$ hyperbolic tetrahedra, if $0 < l_i < \arccosh 3$ for $1\leq i\leq 6$, then the dihedral angle at the longest edge always has a lower bound $\pi/5$ (see Corollary 3.7 in \cite{Feng2022}), and $l$ is a real hyperbolic polyhedral metric (see Theorem 3.9 in \cite{Feng2022}). Using the same discussion, we show that if the value interval is changed to $0 < l_i < \arccosh 2$, the conclusions above still hold. The estimation of lower bound of dihedral angle value is still $\pi/5$, in particular, all dihedral angles are acute angles, i.e. $\phi_i(x)\in (0,1)$. The details are as follows:

Since in the $4$-$0$ case, there is no difference between the status of the six vertices, we just need to prove it for $i=2$ without loss of generality. Set $$\phi(x):=\phi_2(x)=\phi_{12}(x), \quad \alpha(x):=\alpha_2(x)=\alpha_{12}(x).$$
By Theorem 3.6 in \cite{Feng2022}, we have
\begin{eqnarray*}\phi(x)&\leq& k \triangleq \max\{\phi(x_1,2, 2, 1, 2, 2), \phi(x_1, 2, 1, 1, 2, 1),\phi(x_1, 2, 2, 1, 2, 1) \}\\
&=& \max\left\{\frac{-x_1^2 + 8x_1 + 9}{x_1^2 + 8x_1 + 7},
\frac{-x_1^2 + 5x_1 + 5}{x_1^2 + 4x_1 + 4},\frac{-x_1^2 + 6x_1 + 7}{(x_1 + 2)\sqrt{x_1^2 + 8x_1 +7}}\right\}.
\end{eqnarray*}
By the monotonicity of the function, when $x_1 \in (1,2]$ , $k < 1$. Further, when $x_1 = 2$, $k = \frac{7}{9}$. So, \begin{equation}\label{alphalow}
\alpha(x) \geq \arccos {\frac{7}{9}} >\frac{\pi}{5}.
\end{equation}

On the other hand, when $x_i \in (1,2]$,
\begin{eqnarray*}
\phi(x) &\geq & \phi(x_1,x_2,x_3,2,x_5,x_6) \\
&=& \frac{x_2x_3 + x_5x_6 + x_1x_2x_5 + x_1x_3x_6 - 2x ^2_1 + 2}{\sqrt {2x_1x_2x_6 + x^2_1 + x^2_2 + x^2_6 - 1}\sqrt {2x_1x_3x_5 + x^2_1 + x^2_3 + x^2_5 - 1}}\\
& > & \frac{x_2x_3 + x_5x_6 + 2 x_1 - 2x ^2_1 + 2}{\sqrt {2x_1x_2x_6 + x^2_1 + x^2_2 + x^2_6 - 1}\sqrt {2x_1x_3x_5 + x^2_1 + x^2_3 + x^2_5 - 1}}\\
& > & \frac{4 + 2 x_1 - 2x ^2_1}{\sqrt {2x_1x_2x_6 + x^2_1 + x^2_2 + x^2_6 - 1}\sqrt {2x_1x_3x_5 + x^2_1 + x^2_3 + x^2_5 - 1}},
\end{eqnarray*}
Since $x_1 \in (1,2]$,  $4 + 2 x_1 - 2x ^2_1 = -2(x_1 - 2)(x_1 + 1) \geq 0$. Therefore, $\phi(x) > 0$. In a word, when $0 < l_i < \arccosh 2$ , $\phi_i(x)\in (0,1) \subset (-1, 1)$, and $\sigma (l)$ is a hyper-ideal tetrahedron with acute dihedral angles.

\subsection{The induced co-volume function $\widetilde {cov}$ for $3$-$1$ type hyperbolic tetrahedra}
\label{sec-covolume}
In this part, we will introduce a function which is useful in the next section.
Based on the discussion above, let $(\sigma, l)$ be a generalized $3$-$1$ type tetrahedron $\sigma$ with a generalized length $l\in \R^3 \times \R^3_{\ge 0}$. Further assume that the ideal vertex triangular in $(\sigma, l)$ is a Euclidean equilateral triangle with length $1$, where $l_{1i}  = l_{1i} (l_{23} , l_{24}, l_{34} )~(i = 2, 3, 4)$  are defined by the equation (\ref{eq}). Fix $\sigma$ and let $\mathcal{L}^*$ be the space of all such generalized length $l$, then it is diffeomorphic to $\R^3_{> 0}$ by Proposition \ref{inv}.
For each $l\in {\R}^6$, recall $ l_{ij}^+ = \max {\{0, l_{ij}\}}$ by Definition \ref{def-l-l+}. Thus, for $(l_{23} , l_{24}, l_{34}) \in \R^3$,
$$l_{1i} = l_{1i} (l_{23} , l_{24}, l_{34} ) = l_{1i} (l^+_{23} , l^+_{24}, l^+_{34} ), \quad i\in\{2, 3, 4\}$$
is well defined. Thus we can define a new function on $\R^3$ induced by the generalized co-volume function $cov$ as follows:

\begin{defi}[induced co-volume function]
\label{cov-convex}
In a generalized $3$-$1$ type tetrahedron $(\sigma,l)$, the induced co-volume function is defined as
\begin{equation}
\widetilde {cov}(l_{23} , l_{24}, l_{34}) = cov(l_{12} (l_{23} , l_{24}, l_{34} ), l_{13} (l_{23} , l_{24}, l_{34} ),l_{14} (l_{23} , l_{24}, l_{34} ), l_{23} , l_{24}, l_{34} ),
\end{equation}
where $(l_{23} , l_{24}, l_{34}) \in \R^3$.
\end{defi}

\section{The extended combinatorial Ricci flow}
\label{sec-erf}
From now until the end of the article, we assume that $M$ is a compact oriented 3-manifold with boundary, no component of which is a 2-sphere, and $\T$ is an ideal triangulation of $M$. Denote $\partial_t$ by the total boundary components of $M$. For example, $M-\partial_t$ is derived from a volume finite, non-compact, complete hyperbolic 3-manifold with totally geodesic boundary, after forgetting the hyperbolic structure. On the contrary, we want to know if $M-\partial_t$ admits a hyperbolic structure and whether $\T$ is geometric about this hyperbolic structure? The well-known hyperbolization theorem states that $M$ is hyperbolic if and only if it satisfies some topological requirements (such as irreducible, atoroidal and not Seifert fibred). Generally, these topological requirements are not easy to verify. Inspired by the pioneering work of Costantino-Frigerio-Martelli-Petronio in \cite{CFMP}, we show that $M$ is hyperbolic and meanwhile $\T$ is geometric, assuming that the valency on each edge is large enough.

\subsection{Properly glued ideal triangulations}
\label{sec-prop-triangulation}
We refer to \cite{GeJZhang-2024,GeJZhang-2025} for an exact definition of ideal triangulation.
Roughly speaking, an ideal triangulation $\T$ of $M$ is a triangulation with no vertices in $M^{\circ}$. To be precise, let $C(M)$ be the compact $3$-space obtained by coning off each boundary component of $M$ to a point. Suppose $M$ has $s$ torus boundaries and $k$ boundaries of higher genus, then $C(M)$ has $s+k$ cone points $c_1, \cdots, c_s, c_{s+1},\cdots, c_{s+k}$, where the first $s$ cone points correspond to the $s$ torus boundaries.
Note $C(M)-\{c_1,\cdots,c_{s+k}\}$ is homeomorphic to $M^{\circ}$. An ideal triangulation $\T$ of $M$ is a triangulation $\T$ of $C(M)$ such that the vertices of the triangulation are exactly the cone points $c_1,\cdots, c_{s+k}$. By the pioneering work of Moise \cite{Moise}, every compact 3-manifold $M$ can be ideally triangulated. Let $st(c_{s+1},\cdots,c_{s+k})$ be the open star neighbourhood of the vertices $c_{s+1},\cdots, c_{s+k}$ in the second barycentric subdivision of the triangulation $\T$, then $M-\partial_t$ is homeomorphic to $M^{cut}:=C(M)-\{c_1, \cdots, c_s\}-st(c_{s+1},\cdots,c_{s+k})$. For each tetrahedron $\sigma\in \T$, we see $\sigma^{cut}:=\sigma\cap M^{cut}$ is a partially truncated ideal tetrahedron. $\sigma$ is called of type $x$-$(4-x)$ if $\sigma^{cut}$ is truncated at exactly $x$ vertices ($0\leq x\leq 4$), and each such vertex is called a hyper-ideal vertex of $\sigma$. The other $4-x$ vertices are called ideal vertices. Hence each ideal vertex in $\sigma$ corresponds to some cone points in $\{c_1,\cdots,c_s\}$. As a convention, we do not distinguish $\sigma$ and $\sigma^{cut}$. Moreover, each ideal vertex of $\sigma$ is not considered as being contained in $\sigma$ (because it is not contained in $\sigma^{cut}$), and each hyper-ideal vertex of $\sigma$ is represented by the corresponding truncated triangle face in $\sigma^{cut}$.

Let $\{\sigma_1,\cdots, \sigma_t\}$ be the tetrahedra in $\T$ and $\mathscr{T}=\sigma_1 \sqcup\cdots \sqcup \sigma_t$, then $(M,\mathcal{T})$ may be considered as the quotient space $\mathscr{T}/\sim$, where $\sim$ is the gluing relation in $\T$, i.e. a family of affine isomorphisms pairing faces of
tetrahedra in $\mathscr{T}$. The simplexes in $\mathcal{T}$ are equivalent classes of simplexes in $\mathscr{T}$.
Denote $E=E(\mathcal{T})$ and $V=V(\T)$ by the set of edges and vertices in $\mathcal{T}$, which are equivalent classes for vertices and edges of tetrahedra respectively in $\mathscr{T}$ via the gluing. By definition, the \emph{valence} of an edge $e\in E$, denoted by $d_e$, is the number of tetrahedra pasted around $e$, or say, the number of edges in
$\mathscr{T}$ in the equivalent class of $e$.

\begin{figure}[htbp]
	\centering
	\subfloat{\includegraphics[width=.66\columnwidth]{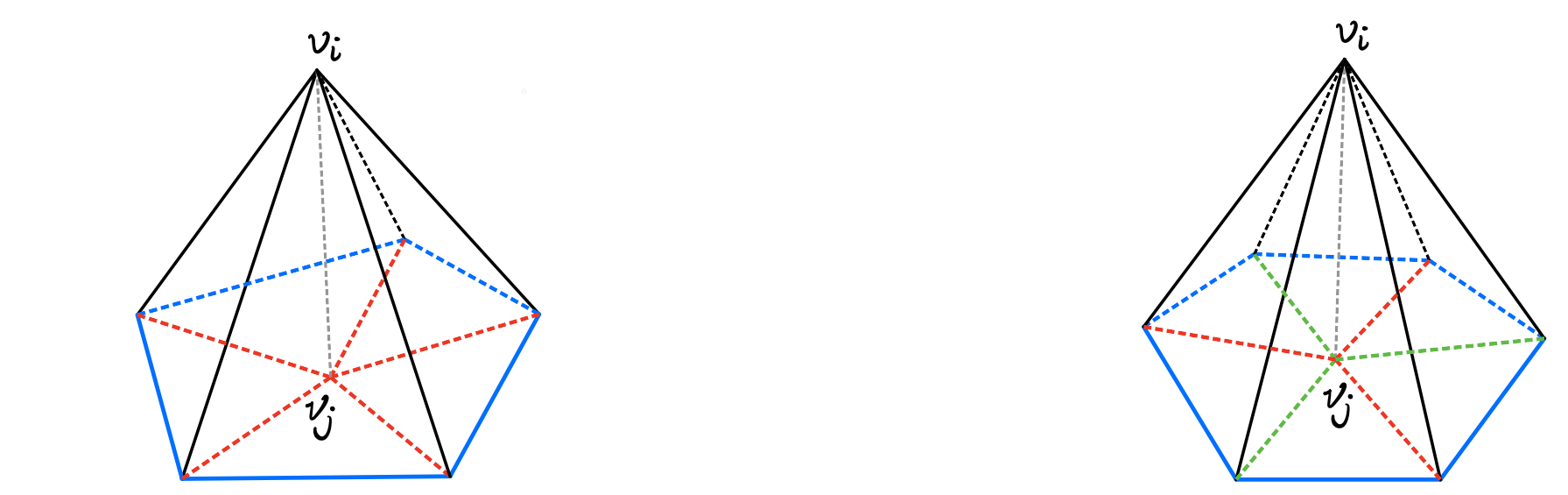}}\hspace{45pt}
	\caption{examples of local proper gluing}
    \label{exam-prop-glue}
\end{figure}

\begin{defi}[properly gluing]
\label{gulingcondition}
The ideal triangulation $\T$ is called \textbf{properly glued}, if each of the tetrahedra in it is of type $3$-$1$ or $4$-$0$, and if two $3$-$1$ tetrahedra $\sigma_{1234}$ and $\sigma_{1'2'3'4'}$ (with ideal vertex $1$ or $1'$ resp.) are glued along the face $123$ and $1'2'3'$ ($k$ and $k'$ are glued for $1\leq k\leq 3$), then the proper gluing requires sides $24$ and $2'4'$ to be glued, $34$ and $3'4'$ to be glued.
 \end{defi}


Recall an \textbf{ideal edge} in a properly glued triangulation $\T$ is a side connecting an ideal vertex and a hyper-ideal vertex. A hyper-ideal edge is a side connecting two hyper-ideal vertices. Moreover, if $\T$ is properly glued at a cusp point $v_i$, there are at most three equivalence classes of hyper-ideal edges corresponding to $v_i$. For example, let $e_{ij}$ be an ideal edge with $v_i$ the ideal vertex, Figure \ref{exam-prop-glue} shows the local structure of properly gluing at $e_{ij}$: there are a total of $10$ ($12$, resp.) hyper-ideal edges in Figure \ref{exam-prop-glue}-left (right, resp.), and the edges of the same color are glued together.


The reason why we consider the proper gluing condition is that when we specify the lengths of the hyper-ideal edges, the mapping $G$ (defined before Proposition \ref{inv}) provides a set of decorated ideal edge lengths. This will be seen in the following sections.

\subsection{ECRF (extended combinatorial Ricci flow)}
\label{sec-ECRF}
\begin{defi}[decorated metric]
\label{def:t1}
A decorated hyperbolic polyhedral metric on $(M, \T)$ is obtained by replacing each tetrahedra by a decorated $3$-$1$ type or $4$-$0$ type hyperbolic tetrahedron according its type and replacing the affine gluing homeomorphisms by isometries preserving the decoration, i.e. gluing decorated hyperbolic tetrahedra along codimension-$1$ faces. By the construction, such a metric is determined by the signed lengths of the edges in $\T$. We denote
$$l=(l(e_1), \dots, l(e_m))$$
by the decorated edge length vector, where $E = \{e_1, \cdots, e_m\}$ is the set of edges in $\T$. In fact, we usually refer to $l$ as a decorated hyperbolic polyhedral metric, or decorated metric for short.
\end{defi}

Denote $\mathcal{L}(M, \mathcal{T})\subset \R^E$ by the set of all decorated metrics on $(M, \mathcal{T})$ parametrized by the decorated edge length vector $l$. Just like each metric tensor determines a curvature tensor in a Riemannian manifold, each decorated metric $l$ determines a combinatorial curvature as follows.

\begin{defi}[curvature]
\label{def:curvature}
Each $l\in\CL(M, \mathcal{T})$ determines a combinatorial Ricci curvature $K\in \mathbb{R}^E$, which is defined as
$$K_e(l)=2\pi-\sum_{e\prec\sigma} \alpha_{\sigma, e}(l)$$
for each edge $e\in E$, where $\alpha_{\sigma, e}$ is the dihedral angle on $e$ in a hyperbolic tetrahedron $\sigma$, and the sum runs over each tetrahedron $\sigma$ with $e$ as one of its edge (denoted by $e\prec \sigma$).
\end{defi}

The following CRF tends to find a decorated metric $l$ with zero curvature, and hence provides a suitable tool to endow a hyperbolic structure on $M$ so that $\T$ is geometric.

\begin{defi}[CRF]
The \textbf{combinatorial Ricci flow} on $(M,\T)$ for decorated metrics $l(t)$ is a system of ODEs with the time parameter $t$ so that
\begin{equation}
\label{eq:luoflow}
\frac{d }{dt}l(t)=K(l(t)),\quad t\geq 0.
\end{equation}
\end{defi}

The solution $l(t)$ may reach $\partial\mathcal{L}(M, \mathcal{T})$ in finite time. To overcome this finite time singularity, we using the extension theory in Section \ref{sec-31type} to extend the definition of $K(l)$ so as it is well-defined even for degenerate $l$, and get the \textbf{ECRF} (extended CRF). By Definition \ref{defi:dihe} and Definition \ref{def-l-l+}, any vector $l \in \mathbb{R}^{E}$ assigns a number to each edge $e\in E$ so that each tetrahedron $\sigma$ becomes a generalized hyperbolic tetrahedron with assigned numbers as decorated edge lengths. Hence we may generalize the definition of a metric.

\begin{defi}[general metric]
\label{gdhpm}
Each $l\in\mathbb{R}^E$ is called a general decorated hyperbolic polyhedral metric (or general metric, in short) on $(M, \mathcal{T})$.
\end{defi}

\begin{defi}[general curvature]
\label{defi:pre}
Each $l\in \mathbb{R}^{E}$ determines a general combinatorial Ricci curvature, which is still denoted by $K\in \mathbb{R}^E$ and is defined as
$$K_e(l)=2\pi-\sum_{e\prec\sigma} \alpha_{\sigma, e}(l^+)$$
for each edge $e\in E$, where $\alpha_{\sigma, e}$ is the extended dihedral angle (see Section \ref{sec-31type}) on $e$ in a hyperbolic tetrahedron $\sigma$, and $l^+$ is defined in Definition \ref{def-l-l+}.
\end{defi}

\begin{rem}
Unlike in Section \ref{sec-introduction}, where we use $\widetilde{K}$ for generalized curvature, from this section onwards, all real and generalized curvature will be represented by
$K$ uniformly for convenience.
\end{rem}

\begin{defi}[ECRF]
The \textbf{ECRF} (extended combinatorial Ricci flow) on $(M,\T)$ is
\begin{equation}
\label{exkl}
\frac{d }{dt}l(t)=K(l(t)),\quad t\geq 0
\end{equation}
or
$$\frac{d}{d t}l_i(t)= K_i(l(t)), \quad i\in E,\; t\geq0.$$
\end{defi}

Formally, (\ref{eq:luoflow}) and (\ref{exkl}) are the same, but they have different meaning. First, the notation $K$ in (\ref{exkl}) represents the generalized curvature. Second, in the ECRF, $l(t)$ is allowed to take values in the whole space $\mathbb{R}^E$, while in the CRF $l(t)$ is only allowed to take values in
$\mathcal{L}(M, \mathcal{T})$.

Suppose $l$ is a decorated metric. For each decorated $3$-$1$ type or $4$-$0$ type hyperbolic tetrahedron $\sigma$, denote $l_{\sigma}\in\mathbb{R}^6$ by the edge length vector of $\sigma$ (i.e. $l$ restricted on the edges of $\sigma$). Let $cov_{\sigma}$ be the generalized co-volume function defined in Corollary \ref{cov}, consider the co-volume functional
\begin{equation}
\label{def-cov}
cov(l) = \sum_{\sigma \in \T} cov_{\sigma}(l_{\sigma}),\quad l\in\mathbb{R}^E
\end{equation}
and
\begin{equation}
\label{def-H}
H(l)=cov(l) - 2\pi\sum_{i=1}^m l_i.
\end{equation}

We list some important results on the ECRF, with details could be found in \cite{Feng2023}.
\begin{prop}\label{H}
The functional $H$ is non-increasing along the ECRF (\ref{exkl}), i.e. for any solution $l(t)$ to the flow (\ref{exkl}),
\[\frac{dH(l(t))}{d t} \le 0.\]
\end{prop}

\tm \label{thm:lt}
For a triangulated compact 3-manifold $(M, \mathcal{T})$ and any initial data $l_0\in\R^E,$ there exists a unique solution $\{l(t)|t\in [0,\infty)\}\subset \R^E$ to the ECRF \eqref{exkl}.

\tmd

\begin{prop}\label{prop:zero}
If a solution $l(t)$ of the ECRF \eqref{exkl} converges to some $\bar l \in \R^E$ as $t \to +\infty$. Then $K(\bar l) = 0.$ In particular, if the solution $l(t)$ converges to $\bar l \in \mathcal {L}(M, \T)$. Then $K(\bar l) = 0$.
\end{prop}

\section{The reduced ECRF}
\label{sec-reduce-ECRF}
Here and after, \textbf{$\T$ always represents an ideal triangulation which is properly glued and every ideal edge has valence 6.} Denote $E=E^*\sqcup E^h$, where $E^*$ is the set of ideal edges, and $E^h$ is the set of hyper-ideal edges. Denote $m^h$ by the number of hyper-ideal edges in $E^h$.

The ECRF on $(M,\T)$ can be rewritten as the following form:
\begin{equation}
\label{eq:fl}
\left\{\begin{array}{cc}
\dfrac{d }{d t}l_i(t)= K_i,\;\; i\in E^h.\\[8pt]
\dfrac{d }{d t}l_j(t)  = K_j, \;\; j\in E^*.
\end{array}\right.
\end{equation}

The presence of ideal vertices and ideal edges makes the ECRF very difficult to handle. We will use the theory developed in subsection \ref{sec-31proper} to introduce a reduced ECRF, which is defined only on $E^h$, and to some extent can be seen as a sub-flow of the original ECRF, but is comparatively easier to deal with because it avoids the ideal edges.
Let $l\subset \mathbb{R}^E$ be a generalized decorated metric, so that each ideal vertex triangle in $\T$ is an Euclidean equilateral triangle with length $1$. Let $\sigma$ be a (maybe degenerated) decorated $3$-$1$ type tetrahedron. Since the ideal vertex triangle in $(\sigma, l_{\sigma})$ is equilateral, the dihedral angle on each ideal edge is always $\pi/6$. By the properly glued assumption on $\T$ that each ideal edge has valence $6$, the generalized curvature is always $0$ at each ideal edge. This shows that each ideal edge already has zero curvature (i.e., hyperbolic metric), and we don't need to evolve the $l^*$ part with the ECRF (\ref{eq:fl}). We just need to consider the $l^h$ part in the ECRF (\ref{eq:fl}).

\begin{lem}
\label{lem-l^h-determine-l}
For each $l^h\in\mathbb{R}^{E^h}_{>0}$, $(M,\T)$ admits a unique generalized metric $l$, so that $l^h$ is the hyper-ideal part of $l$, that is, the projecting of $l\in \mathbb{R}^E$ onto $\mathbb{R}^{E^h}$, and each ideal vertex triangle is equilateral with length 1.
\end{lem}
\begin{proof}
Take a $3$-$1$ tetrahedron $\sigma\in\T$ and an ideal edge $e\prec\sigma$, by Proposition \ref{surjection} and Corollary \ref{prop-equa-length-1}, the decorated edge length $l(e)$ is exist and uniquely determined by $l^h_{\sigma}=(l^h_{\sigma,1},l^h_{\sigma,3},l^h_{\sigma,5})$ (i.e. the edge lengthes of the three hyper-ideal edges in $\sigma$).
In addition, if $e$ is a common ideal edge of $\sigma$ and $\sigma'$, the properly glued assumption on $\T$ ensures that $l^h_{\sigma}=l^h_{\sigma'}$. So there is no difference whether $l(e)$ is determined in $\sigma$ or $\sigma'$. So $l=(l^*, l^h)$ is uniquely determined by its $l^h$ part.
\end{proof}

By Lemma \ref{lem-l^h-determine-l}, for each $l^h=\{l_i\}_{i\in E^h}\in\mathbb{R}^{E^h}_{>0}$, we can use the formula (\ref{eq}) to calculate $l^*$, and thus reconstruct the generalized metric $l$. This shows that the functions $K$, $cov$ and $H$ could be viewed as functions defined for $l^h\in\mathbb{R}^{E^h}_{>0}$.

\begin{defi}[the reduced ECRF]
The \textbf{reduced ECRF} on $(M,\T)$ is defined as
\begin{equation}
\label{equation-reduced-ECRF}
\frac{d}{dt}l_i(t)=K_i(l^h)=K_i(l^*, l^h),\;\; i\in E^h, \;\; l^h\in\mathbb{R}^{E^h}_{>0}.
\end{equation}
\end{defi}

Recall the functions $cov(l)$ and $H(l)$ are defined in (\ref{def-cov} and \ref{def-H}) respectively. For the sake of subsequent elaboration, we rewrite the functions $cov$ and $H$ as
\begin{equation}
{cov}^h(l^h)=cov(l^*,l^h)=cov(l) = \sum_{\sigma\in\T_{3,1}} \widetilde {cov}_{\sigma}(l^h_{\sigma}) + \sum_{\sigma\in\T_{4,0}}cov_{\sigma}(l_{\sigma})
\end{equation}
and
\begin{equation}
H^h(l^h)=H(l^*,l^h)=H(l)= cov(l) - 2\pi\sum_{i\in E} l_i,
\end{equation}
where $\T_{3,1}$ ($\T_{4,0}$, resp.) represents the set of $3$-$1$ type ($4$-$0$ type, resp.) tetrahedra in $\T$.

\begin{lem}
For $(M,\T)$, let $l^h\in \mathcal{L}^h(M,\T)$ be a hyper-ideal edge length vector, then
\[\frac{\partial H^h}{\partial l_i} = - K_i, \quad e_i\in E^h.\]
\end{lem}
\begin{proof}
For any $e_i\in E^h$, if the hyper-ideal edge $e_i$ is not a side of any $3$-$1$ type tetrahedron in $\T$, then (4-1) in \cite{Feng2022} had proved the conclusion.
On the contrary, $e_i$ is a side of some $3$-$1$ type hyperbolic tetrahedra, take such a $3$-$1$ tetrahedron $(\sigma, l_{\sigma})$, with
$$l_{\sigma}=(l_{12}(l_{23}, l_{24},  l_{34}), l_{13}(l_{23}, l_{24},  l_{34}), l_{14}(l_{23}, l_{24},  l_{34}), l_{23}, l_{24},  l_{34}),$$
where $(l_{23}, l_{24},  l_{34})=l^h_{\sigma} \in \R^3$, then
\begin{eqnarray}\label{pcov-23}
\frac{\partial \widetilde {cov}_{\sigma}(l^h_{\sigma})}{\partial l_{23}}
&=& \frac{\partial cov(l_{\sigma})}{\partial l_{23}} \nonumber\\
&=& \frac{\partial cov}{\partial l_{23}} + \frac{\partial cov}{\partial l_{12}}\frac{\partial l_{12}}{\partial l_{23}} + \frac{\partial cov}{\partial l_{13}}\frac{\partial l_{13}}{\partial l_{23}}+ \frac{\partial cov}{\partial l_{14}}\frac{\partial l_{14}}{\partial l_{23}}\\
&=& \alpha_{23} + \alpha_{12}\frac{\partial l_{12}}{\partial l_{23}} + \alpha_{13}\frac{\partial l_{13}}{\partial l_{23}} + \alpha_{14}\frac{\partial l_{14}}{\partial l_{23}}  \label{eq:newderiv} \nonumber\\
&=& \alpha_{23} + \frac{\pi}{3}\frac{\partial}{\partial l_{23}}\big(l_{12}+l_{13}+l_{14}\big). \nonumber
\end{eqnarray}

Similarly, we also have

\begin{equation}\label{pcov-24}
\frac{\partial \widetilde {cov}_{\sigma}(l^h_{\sigma})}{\partial l_{24}} =  \alpha_{24} +\frac{\pi}{3}\frac{\partial}{\partial l_{24}}\big(l_{12}+l_{13}+l_{14}\big).
\end{equation}
and
\begin{equation}\label{pcov-34}
\frac{\partial \widetilde {cov}_{\sigma}(l^h_{\sigma})}{\partial l_{34}} =  \alpha_{34} + \frac{\pi}{3}\frac{\partial}{\partial l_{34}}\big(l_{12}+l_{13}+l_{14}\big).
\end{equation}

Substitute (\ref{pcov-23}), (\ref{pcov-24}) and (\ref{pcov-34}) into the following calculation of $\partial H^h/\partial l_i$, we have

\begin{eqnarray}
\frac{\partial H^h}{\partial l_i}
&=& \frac{\partial {cov}^h(l^h)}{\partial l_i} - 2\pi\sum_{j\in E}\frac{\partial l_j}{\partial l_i}\nonumber\\
&=&\frac{\partial }{\partial l_i}\sum_{\sigma\in\T_{3,1}}\widetilde{cov}_{\sigma}(l^h_{\sigma})+\frac{\partial}{\partial l_i}\sum_{\sigma\in\T_{4,0}}cov_{\sigma}(l_{\sigma})-2\pi(1+ \sum_{j\in E^*}\frac{\partial l_j}{\partial l_i})\nonumber\\
&=&\sum_{\sigma\in\T_{3,1}}\alpha(\sigma, e_i)+\frac{\pi}{3}\sum_{\sigma\in\T_{3,1}}\frac{\partial}{\partial l_i}\big(l_{\sigma,12}+l_{\sigma,13}+l_{\sigma,14}\big)+\sum_{\sigma\in\T_{4,0}}\alpha(\sigma, e_i)-2\pi-2\pi\sum_{j\in E^*}\frac{\partial l_j}{\partial l_i}\nonumber\\
&=&\bigg(\sum_{\sigma\in\T}\alpha(\sigma, e_i)-2\pi\bigg)+\bigg(\frac{\pi}{3}\sum_{\sigma\in\T_{3,1}}\frac{\partial}{\partial l_i}\big(l_{\sigma,12}+l_{\sigma,13}+l_{\sigma,14}\big)-2\pi\sum_{j\in E^*}\frac{\partial l_j}{\partial l_i}\bigg)\nonumber\\
&=& \sum_{\sigma\in\T}\alpha(\sigma, e_i) - 2\pi+0\nonumber\\
&=& -{K}_i,\nonumber
\end{eqnarray}
where $\alpha(\sigma, e_i)$ is the dihedral angle at $e_i$ in $\sigma$ (it is zero if $e_i$ is not a side of $\sigma$), and we used the fact that each ideal edge $j\in E^*$ has zero curvature to derive the penultimate equal sign.
\end{proof}

Consequently, the reduced ECRF is a negative gradient flow, that is,
\begin{equation}
\frac{d l^h(t)}{dt}=K(l^h)=-\nabla_{l^h}H^h,
\end{equation}
which implies the following
\begin{prop}\label{H-neg}
The functional $H^h$ is non-increasing along the reduced ECRF (\ref{equation-reduced-ECRF}), that is, for any solution $l^h(t)$ to the reduced ECRF (\ref{equation-reduced-ECRF}),
\[\frac{d H^h(l^h(t))}{d t}=-\|K\|^2\le 0.\]
\end{prop}

Consider the reduced ECRF (\ref{equation-reduced-ECRF}) with an initial metric $l^h_0 \in \R^{E^h}_{> 0}$, we have

\tm\label{thm:lt1}
For any initial data $l^h_0 \in \R^{E^h}_{> 0}$, there exists a solution $\{l^h(t)|t\in [0,\infty)\}\subset \R^{E^h}_{>0}$ to the reduced ECRF (\ref{equation-reduced-ECRF}).
\tmd
\begin{proof}
Since $M$ is compact, the triangulation $\T$ is finite, hence all curvatures $K_i$ are uniformly bounded. Hence $\|l_i(t)\|=\|\int^t K_i\|\leq \exp\{c_1t+c_2\}$, implying that $l^h(t)$ will not reach the infinity boundary of $\R^{E^h}_{> 0}$. In addition, by Proposition \ref{prop:c0} in this article and Proposition 3.10 in \cite{Feng2022}, if some $l_i(t)>0$ is close to $0$, then $\alpha(\sigma, e_i)<\epsilon$ uniformly for all $\sigma\succ e_i$. It follows
$$K_i=2\pi-\sum\alpha(\sigma, e_i)>2\pi-\epsilon d_{\rm{max}}>0,$$
where $d_{\rm{max}}$ is the maximal edge valence in $\T$, and $\epsilon=2\pi/d_{\rm{max}}$. Then $\frac{dl_i}{dt}=K_i>0$, so $l_i(t)$ is increasing whenever it reaches zero, implying that $l(t)$ will not reach the zero boundary of $\R^{E^h}_{> 0}$ along the reduced ECRF. Hence the solution $l^h(t)$ to (\ref{equation-reduced-ECRF}) exists for all time $t\geq0$.
\end{proof}

\tm
\label{thm:expc1}
If a solution $l^h(t)$ of the reduced ECRF (\ref{equation-reduced-ECRF}) converges to some $\overline{l^h} \in \R^{E^h}_{>0}$ as $t \to +\infty$. Then $K(\overline{l^h}) = 0.$
\tmd
Since the proof of Theorem \ref{thm:expc1} is almost the same with the proof of Theorem 4.9 in \cite{Feng2022}), we omit the details.

\section{Proof of Main Theorem \ref{mainthm}}
\label{sec-proofthm}
We will prove Theorem \ref{mainthm} in this section, and $(M, \mathcal{T})$ satisfies the assumptions in Section \ref{sec-reduce-ECRF}.

Although $H$ is a convex function of $l$, we do not know whether $H^h$ is convex as a function of $l^h$, which seems a bit subtle. So, our proof strategy in \cite{Feng2022} and the convex optimization methods used in \cite{Feng2023} can not be generalized here directly. Hence a convergent $l^h(t)$ is hard to obtain generally. However, a subsequence convergence of $l(t)$ can be obtained under suitable assumption on $\T$, and this is sufficient for finding hyperbolic structures and geometric triangulations.

Note that the function $H^h(l^h(t))$ is non-increasing, we may consider getting a uniform bound for the solution $l^h(t)$ to (\ref{equation-reduced-ECRF}).
We give the upper bound estimation firstly.

\begin{theo}\label{LONG}
Suppose the valence of each hyper-ideal edge is at least 11,  then for any initial metric
$l^h_0 \in (0,\arccosh 2)^{E^h}$, the solution $\{l^h(t)|t\in[0,\infty)\}$ to the reduced ECRF (\ref{equation-reduced-ECRF}) satisfies
$$l^h(t)\in (0,\arccosh 2)^{E^h} \subset \R^{E^h}_{>0}, \quad \forall t\geq0.$$
\end{theo}

\begin{proof}
We need to prove $l_i(t)<\arccosh 2$ for all $t\geq 0 $ and all $i\in E^h$. If it is not true, then there is a non-empty closed set in $\R$,
$$S:=\{t\in [0,\infty): \max_{i\in E^h} l_i(t)\geq \arccosh 2\}.$$

Set $t_0:=\inf S<\infty$. Obviously $0\not\in S$, then $t_0>0$, and there exists a $e_j\in E^h$ so that
$$l_j(t_0)=\max_{i\in E^h} l_i(t_0).$$
The equation (\ref{equation-reduced-ECRF}) implies that $l_j(t)$ is $C^1$-smooth.
Note that $l_j(t_0)=\arccosh 2$ and $$l_j(t)\leq \max_{i\in E^h} l_i(t)< \arccosh 2,\quad \forall\ t<t_0,$$ it follows
$$l_j'(t_0)\geq 0.$$

By Corollary ~\ref{lem:longest} and (\ref{alphalow}), since $l_j(t_0)=\max_{i\in E^h} l_i(t_0)=\arccosh 2$,  for any $\sigma\succ e_j$,
$$\alpha(\sigma,e_j) > \frac{2\pi}{11}. $$

Using $d_j \ge 11$,  we see
$$ K_j(l^h(t_0))=2\pi-\sum_{\sigma\succ e_j}\alpha(\sigma, e_j)<2\pi-\frac{2\pi}{11}d_j\leq 0. $$
This contradicts with $0\leq l_j'(t_0)=K_j (l^h(t_0)) <0$ , This proves the statement.
\end{proof}

\begin{rem}
In fact, since the dihedral angle on a hyper-ideal tetrahedron is at least greater than $\frac{\pi}{5}$ when the longest edge of the hyper-ideal tetrahedron is $\arccosh 2$, as discussed in the proof, if an edge $e$ only belongs to the hyper-ideal tetrahedron, $d_e \geq 10$ is enough.
\end{rem}

Next, we give the lower bound estimate:

\begin{theo}\label{sbd}
Let $\{l^h(t)|t\in[0,\infty)\}\subset \R^{E^h}_{> 0}$ be a solution to the reduced ECRF (\ref{equation-reduced-ECRF}) with the initial metric $l^h_0 \in (0,\arccosh 2)^{E^h}$. Suppose that there exists a constant $C>0$ such that $$l^h(t)\in (0,C]^{E^h},\quad \forall t\geq 0.$$ Then there exists a positive constant $c(l_0,C,d_{\rm{max}}),$ depending on $l_0,$ $C$ and $d_{\rm{max}}$ such that
\begin{equation}\label{eq:spq1}l(t)\in (c,C]^{E^h},\quad \forall t\geq 0.\end{equation}
\end{theo}

\begin{proof}
Set $\epsilon_0:=\pi/d_{\rm{max}}$, and
\[c:=\min\Big\{\frac12 \min_{i\in {E^h}}l_{0,i},\arccosh(1+\delta)\Big\},\]
where $\delta=\delta(\cosh C,\epsilon_0)$  is the constant given in Proposition~\ref{prop:c0}. We will show \eqref{eq:spq1} holds for the above constant $c$.

Suppose that it is not true, then
$$Q:=\{t\in [0,\infty): \min_{i\in {E^h}} l_i(t)\leq c\}$$ is a non-empty closed set in $\R$. Set $t_0:=\inf Q.$ Since \eqref{eq:spq1} holds for $t=0$,  $0 \not\in Q.$ Hence $0<t_0<\infty,$  There exists some $j\in E^h$ such that
$$l_j(t_0)=\min_{i\in E^h} l_i(t_0).$$

Note that $l_j(t_0)=c$ and $$l_j(t)\geq \min_{i\in E^h} l_i(t)>c,\quad \forall\ t < t_0.$$ This yields that $$l_j'(t_0)\leq 0.$$

Note that $l(t_0)\in (0,C]^{E^h}$ and $$\cosh(l_j(t_0))=\cosh  c \le 1+\delta. $$
Let $e$ be a hyper-ideal edge in a $3$-$1$ type hyperbolic tetrahedron realizing the length $l_j(t_0)$, by Proposition \ref{prop:c0}, for all $3$-$1$ type tetrahedra $\sigma\succ e$, we have
$$\alpha(e)<\epsilon_0.$$
Considering this is also true for $4$-$0$ type tetrahedra by Proposition 3.10 in \cite{Feng2022}, thus,
$$K_j(l(t_0))=2\pi-\sum_{\sigma\succ e}\alpha(\sigma,e)>2\pi-d_j\epsilon_0\geq \pi> 0.$$
This yields the following contradiction,
$$0\geq l_j'(t_0)=K_j (l(t_0)) > 0.$$

This proves the result.
\end{proof}
By the same argument as above, we can prove a quantitative lower bound estimate for the
solution to the reduced ECRF (\ref{equation-reduced-ECRF}). We denote $\mathbbm{1}_{E^h}$ by the constant function $1$ on $E^h$, i.e.
$$\mathbbm{1}_{E^h}=(1,\cdots,1).$$

\begin{theo}\label{bd}
Suppose the valence of each hyper-ideal edge is at least 11, and $\{l(t)|t\in[0,\infty)\}\subset \R^{E^h}_{> 0}$ is a solution to the reduced ECRF (\ref{equation-reduced-ECRF})  with an initial metric $l^h_0 \in (0,\arccosh 2)^{E^h}$. Then there exists a positive constant $c(l_0,d_{\rm{max}})$ such that
\begin{equation}
\label{eq:pq1}
l^h(t)\in (c,\arccosh 2)^{E^h}, \quad \forall t\geq 0.
\end{equation}
In particular, if $l_0=\frac{1}{2}\arccosh 2\cdot\mathbbm{1}_{E^h},$ then $c$ can be chosen as $c(l_0,d_{\rm{max}})=1/3 d_{\rm{max}}$.
\end{theo}
\begin{proof}

By Theorem~\ref{LONG}, $$l^h(t)\in (0,\arccosh 2)^{E^h},\quad \forall t\geq 0.$$
Applying Theorem~\ref{sbd} for $C=\arccosh 2,$ we get the result \eqref{eq:pq1}.

Let $l_0=\frac{1}{2}\arccosh 2\cdot\mathbbm{1}_{E^h}.$ Set $\epsilon_0:=\pi/d_{\rm{max}}.$
By the same argument as in the proof of Theorem~\ref{sbd}, we may prove that $$l(t)\in (c_1,\arccosh 2)^{E^h}, \quad \forall t\geq 0,$$ where $c_1=\min\left\{\frac{1}{4}\arccosh 2,\arccosh(1+\delta)\right\}.$ Here $\delta=\delta(2,\epsilon_0)=\frac{1}{5d_{\rm{max}}^2}\leq \frac{1}{500}.$
Note that $$c_1=\arccosh(1+\delta)\geq \sqrt{\delta}\geq \frac{1}{3d_{\rm{max}}}.$$
Hence we can choose $c=\frac{1}{3d_{\rm{max}}}$ as a new lower bound in \eqref{eq:pq1}.
\end{proof}

Now we are ready to prove Theorem~\ref{mainthm}.
\begin{proof}[Proof of Theorem~\ref{mainthm}]
For any initial data $l^h_0\in (0,\arccosh 2)^{E^h},$ let $l^h(t)$ be the solution to (\ref{equation-reduced-ECRF}). By Theorem~\ref{LONG} and Theorem~\ref{bd}, there exist a constant $C_1=c(l_0,d_{\rm{max}})$ such that
\begin{equation}\label{eq:mm1}l(t)\in (C_1,C_2)^{E^h},\quad \forall t\geq 0,\end{equation} where $C_2=\arccosh 2.$

By \eqref{eq:mm1} and the continuity of $H^h$, the set
$\{H^h(l^h(t)) : t \ge 0\}\subset \R$ is bounded.
By the monotonicity of $H^h(l(t)),$ the following limit exists and is finite,
\[\lim_{t \to \infty} H^h(l(t)) = C.\]
Consider the sequence $\{H^h(l^h(n))\}^\infty_{n=1}$. By the mean value theorem, for any
$n \ge 1$, there exists a $t_n \in (n, n + 1)$ such that
\beq\label{eq:aaa1}
H^h(l^h(n + 1)) - H^h(l^h(n)) =\frac{d}{dt}\big |_{t = t_n} H^h(l^h(t)) = - \sum_i K_i^2(l^h(t_n)).
\eeq
Note that $\lim_{n \to \infty} H^h(l^h(n + 1)) - H^h(l^h(n)) = 0$ and we have the estimate \eqref{eq:mm1}. Passing to the limit, $n\to\infty,$ in \eqref{eq:aaa1}, we have
\begin{equation}\label{eq:mm2}\lim_{n\to \infty}|K_i(l^h(t_n))|=0,\quad \forall i\in {E^h}.\end{equation} By \eqref{eq:mm1}, there exist a subsequence of $\{t_n\}_{n=1}^\infty,$ denoted by $\{t_{n_k}\}_{k=1}^\infty,$ and $l^h_\infty\in [C_1,C_2]^{E^h}$ such that
\begin{equation}\label{eq:mm4}l^h(t_{n_k})\to l^h_\infty,\quad k\to\infty.\end{equation}
By the continuity of $K(l^h)$ and \eqref{eq:mm2}, we have
$$K(l^h_\infty)=0.$$

To obtain a refined estimate for $l^h_\infty,$ we choose a specific initial metric $\overline{l^h_0}=\frac{1}{2}\arccosh 2\cdot\mathbbm{1}_{E^h}$.
Let $\overline{l^h}(t)$ be a solution to the reduced ECRF (\ref{equation-reduced-ECRF}). By Theorem~\ref{LONG} and Theorem~\ref{bd}, we have
 $$\overline{l^h}(t)\in ((3 d_{\rm{max}})^{-1}, \arccosh 2)^{E^h},\quad \forall\ t\geq 0.$$ Passing to the limit, we get
 $$l^h_\infty\in [(3 d_{\rm{max}})^{-1}, \arccosh 2]^{E^h}.$$

Now, since $\T$ satisfies the proper gluing condition, we can get a decorated hyperbolic polyhedron metric $\bar l_{\infty}$ such that its every hyper-ideal edge component is exactly the corresponding component in $l^h_{\infty}$ and its every ideal edge component is calculated from (\ref{eq}). Then, $\bar l_{\infty}$, as the generalized decorated hyperbolic polyhedron metric, satisfies that every ideal vertex triangle is regular triangle with length $1$ and satisfies the curvature on each edge is zero. Further, by Theorem ~\ref{Range} and Proposition ~\ref{realmetric}, $\bar l_{\infty}$ is a realizable decorated hyperbolic metric.
At last, based on the rigidity of decorated hyperbolic polyhedral metric proved by the first two authors, see Theorem 1.1 in \cite{Feng2023}, this zero curvature metric $\bar l_{\infty}$ is unique.
\end{proof}

\section{Appendix}

Now we give an example of a manifold $(M, \T)$ that satisfies the conditions of Theorem \ref{mainthm}. We can compute $\mathcal{T}$ with only three vertices, where the links at the vertex $v_1$ and the vertex $v_3$ are tori, and the link at the vertex $v_2$ is a surface with genus six.
The degree of the edge connected $v_1$ to $v_2$ is 6 , the degree of the edge connected $v_2$ to $v_3$ is 6 and the degree of the edge connected $v_2$ to $v_2$ is 36.
The following table provides a triangulation.
In each tetrahedron, we denote its vertices by $0, 1, 2, 3$, where $0$ represents the ideal vertex and the rest represents the hyper-ideal vertex. The first column of the table lists twelve tetrahedra, labeled 0 to 11, and the remaining elements in the table indicate how the faces are glued together.

   $$ \begin{array}{|c|c|c|c|c|}\hline \text { Tetrahedron } & \text { Face } 012 & \text { Face } 013 & \text { Face } 023 & \text { Face } 123 \\\hline 0 & 3(012) & 5(013) & 1 (023)  & 2(213) \\\hline 1 & 4(012) & 2(013) & 0(023) & 11  (123)  \\\hline 2 & 5(012) & 1(013) & 3(023) & 0(213) \\\hline 3 & 0(012) & 4(013) & 2(023) & 4 (231)\\\hline 4 & 1  (012)  & 3(013) & 5(023) & 3  (312)  \\\hline 5 & 2(012) & 0(013) & 4(023) & 7   (123)  \\\hline 6 & 9(012) & 11  (013)  & 7(023) &8(213) \\\hline 7 & 10 (012) & 8(013) & 6(023) & 5(123) \\\hline 8 & 11  (012)  & 7   (013)  & 9(023) & 6(213) \\\hline 9 & 6(012) & 10(013) & 8(023) & 10(231) \\\hline 10 & 7  (012)  & 9(013) & 11(023) & 9  (312)  \\\hline 11 & 8(012) & 6(013) & 10(023) & 1  (123) \\\hline
    \end{array}$$

\noindent Ke Feng, kefeng@uestc.edu.cn\\
\emph{School of Mathematical Sciences, University of Electronic Science and Technology of China, Sichuan 611731, P. R. China}\\[2pt]

\noindent Huabin Ge, hbge@ruc.edu.cn\\
\emph{School of Mathematics, Renmin University of China, Beijing 100872, P. R. China}\\[2pt]

\noindent Yunpeng Meng, 2230501005@cnu.edu.cn\\
\emph{School of Mathematical Sciences, Capital Normal University, Beijing 100048, P. R. China}


\begin{thebibliography}{99}

\bibitem{Aki}
H. Akiyoshi, \emph{On the Ford domains of once-punctured torus groups}, from: ``Hyperbolic spaces and related topics (Japanese) (Kyoto, 1998)", S\={u}rikaisekikenky\={u}sho K\={o}ky\={u}roku 1104 (1999), 109-121.

\bibitem{Aki-SWY}
H. Akiyoshi, M. Sakuma, M. Wada, Y. Yamashita, \emph{J{\o}rgensen's picture of punctured torus groups and its refinement}, from: ``Kleinian groups and hyperbolic 3-manifolds (Warwick, 2001)", (Y. Komori, V. Markovic, C. Series, editors), London Math. Soc. Lecture Note Ser. 299, Cambridge Univ. Press (2003), 247-273.

\bibitem{Aki-SWY-2}
H. Akiyoshi, M. Sakuma, M. Wada, Y. Yamashita, \emph{Punctured torus groups and 2-bridge knot groups. I}, Lecture Notes in Math. 1909, Springer, Berlin (2007).


\bibitem{Bobenko2015}
  A. I. Bobenko, U. Pinkall, B. A. Springborn, \emph{Discrete conformal maps and ideal hyperbolic polyhedra}, Geom. Topol. 19 (2015), no. 4, 2155-2215.

\bibitem{Chow-Luo}
  B. Chow, F. Luo, \emph{Combinatorial Ricci flows on surfaces}, J. Differential Geometry, 63 (2003), 97-129.


  \bibitem{CFMP}
  F. Costantino, R. Frigerio, B. Martelli, C. Petronio, \emph{Triangulations of 3-manifolds, hyperbolic relative handlebodies, and Dehn filling}, Comment. Math. Helv. 82 (2007), 903-933.

  \bibitem{CDGW2016}
  M. Culler, N. M. Dunfield, M. Goerner, J. R. Weeks, \emph{SnapPy, a computer program for studying the geometry and topology of 3-manifolds}, Available at http://snappy.computop.org, 2016.

  \bibitem{Feng2023}
  K. Feng, H. Ge, \emph{Combinatorial Ricci flow and Thurston's triangulation conjecture}, arXiv:2502.06497, 2025.


  \bibitem{Feng2022}
  K. Feng, H. Ge, B. Hua, \emph{Combinatorial Ricci flows and the hyperbolization of a class of compact 3-manifolds}, Geom. Topol. 26 (2022), no. 3, 1349-1384.

  \bibitem{Feng2022-1}
  K. Feng, H. Ge, B. Hua, X. Xu, \emph{Combinatorial Ricci flows with applications to the hyperbolization of cusped 3-manifolds},  Int. Math. Res. Not. (2022), no. 20, 15549-15573.

  \bibitem{FengGLiu}
  K. Feng, H. Ge, C. Liu, \emph{Rigidity, volume and angle structures of 1-3 type hyperbolic polyhedral 3-manifolds}, arXiv:2501.08081, 2025.

  \bibitem{GeHua}
  H. Ge, B. Hua, \emph{3-dimensional combinatorial Yamabe flow in hyperbolic background geometry}, Trans. Amer. Math. Soc. 373 (2020), no. 7, 5111-5140.

  \bibitem{GeJZhang-2024}
  H. Ge, L. Jia, F. Zhang, \emph{Angle structure on general hyperbolic 3-manifolds}, arXiv:2408.14003, 2024.

  \bibitem{GeJZhang-2025}
  H. Ge, L. Jia, F. Zhang, \emph{Angle structures on pseudo 3-manifolds}, arXiv:2502.11397, 2025.

  \bibitem{GeJShen}

  H. Ge, W. Jiang, L. Shen, \emph{On the deformation of ball packings}, Adv. Math. 398 (2022), Paper No. 108192, 44 pp.

  \bibitem{GeLin}
  H. Ge, A. Lin, \emph{The character of Thurston's circle packings}, Sci. China Math. 67 (2024), no. 7, 1623-1640.


  \bibitem{Gue-1}
  F. Gu\'{e}ritaud, \emph{On canonical triangulations of once-punctured torus bundles and two-bridge link complements}, Geom. Topol. 10 (2006), 1239-1284, With an appendix by D. Futer.

  \bibitem{Gue-2}
  F. Gu\'{e}ritaud, \emph{Triangulated cores of punctured-torus groups}, J. Differential Geom. 81 (2009), 91-142.

  \bibitem{GS2010} F. Gueritaud, S. Schleimer, \emph{Canonical triangulations of Dehn fillings}, Geom. Topol. 14 (2010), no. 1, 193-242.

  \bibitem{Ham-P} S. L. Ham, J. S. Purcell, \emph{Geometric triangulations and highly twisted links}, Algebr. Geom. Topol. 23 (2023), 1399-1462.


 \bibitem{Jorgen}
  T. J{\o}rgensen, \emph{On pairs of once-punctured tori}, from: ``Kleinian groups and hyperbolic 3-manifolds (Warwick, 2001)", (Y. Komori, V. Markovic, C. Series, editors), London Math. Soc. Lecture Note Ser. 299, Cambridge Univ. Press (2003), 183-207.


  \bibitem{Lackenby2000}
  M. Lackenby, \emph{Word hyperbolic Dehn surgery}, Invent. Math. 140 (2000), no. 2, 243-282.
  
   \bibitem{Lackenby-AGT}
  M. Lackenby, \emph{Analgorithm to determine the Heegaard genus of simple 3-manifolds with non-empty boundary}, Alg. Geom. Top. 8 (2008), 911-934.

    \bibitem{Lackenby}
  M. Lackenby, \emph{The canonical decomposition of once-punctured torus bundles}, Comment. Math. Helv. 78 (2003), 363-384.


  \bibitem{Luo2005}
  F. Luo, \emph{A combinatorial curvature flow for compact 3-manifolds with boundary}, Electron. Res. Announc. Amer. Math. Soc. 11 (2005), 12-20.


  \bibitem{Luo2011}
  F. Luo, \emph{Rigidity of polyhedral surfaces, III}, Geom. Topol. 15 (2011), no. 4, 2299-2319.


  \bibitem{Luo2018}
  F. Luo, T. Yang, \emph{Volume and rigidity of hyperbolic polyhedral 3-manifolds}, J. Topol. 11 (2018), no. 1, 1-29.

  \bibitem{Moise}
   E. E. Moise, \emph{Affine structures in 3-manifolds, V: The triangulation theorem and Hauptvermutung}, Ann. of Math. 56 (1952) 96-114.

 
 \bibitem{Nimer}
  B. Nimershiem, \emph{Geometric triangulations of a family of hyperbolic 3-braids}, Algebr. Geom. Topol. 23 (2023), no.9, 4309-4348.


  \bibitem{Rivin1994}
  I. Rivin, \emph{Euclidean structures on simplicial surfaces and hyperbolic volume}, Ann. of Math. (2) 139 (1994), no. 3, 553-580.

  \bibitem{Rivin2003}
  I. Rivin, \emph{Combinatorial optimization in geometry}, Adv. in Appl. Math. 31 (2003), no. 1, 242-271.

  \bibitem{Thurston2022}
  W. P. Thurston, \emph{The geometry and topology of three-manifolds}, lecture notes, Princeton Univ. (1980) Available at http://library.msri.org/books/gt3m/.\\[8pt]

\end{thebibliography}
\end{document}